%% file: CiJa2x-StokesGrad-hal.tex
\documentclass[preprint,12pt]{elsarticle}
\usepackage[utf8]{inputenc}   
\usepackage[T1]{fontenc}
\usepackage[french,english]{babel}
\usepackage{aeguill}
\usepackage{amsmath}%
\usepackage{amsthm}%
\usepackage{amsfonts}%
\usepackage{amssymb}%
\usepackage{graphicx}
\usepackage{color}
\usepackage{dsfont}
\usepackage{bbm}
\usepackage{pgfplots}
\usepackage{tikz}
\usetikzlibrary{calc,intersections}
\usepackage{url}
\usepackage[normalem]{ulem}
\usepackage{multirow}
\usepackage{algorithm}
\usepackage{algorithmic}
\newcommand{\im}{\operatorname{im}}
\newcommand{\eps}{\varepsilon}
\newcommand{\Ccal}{\mathcal{C}}

\newcommand{\Ical}{\mathcal{I}}
\newcommand{\Jcal}{\mathcal{J}}

\newcommand{\Lcal}{\mathcal{L}}

\newcommand{\Tcal}{\mathcal{T}}
\newcommand{\Xcal}{\mathcal{X}}
\newcommand{\bbA}{\mathbb{A}}
\newcommand{\bbB}{\mathbb{B}}
\newcommand{\R}{\mathbb{R}}

\newcommand{\bbK}{\mathbb{K}}
\newcommand{\bbL}{\mathbb{L}}
\newcommand{\bbM}{\mathbb{M}}
\newcommand{\N}{\mathbb{N}}
\newcommand{\bL}{\mathbf{L}}

\newcommand{\bH}{\mathbf{H}}
\newcommand{\bP}{\mathbf{P}}
\newcommand{\bV}{\mathbf{V}}

\newcommand{\evec}{\mathbf{e}}
\newcommand{\fvec}{\mathbf{f}}
\newcommand{\gvec}{\mathbf{g}}
\newcommand{\nvec}{\mathbf{n}}

\newcommand{\uvec}{\mathbf{u}}
\newcommand{\vvec}{\mathbf{v}}
\newcommand{\wvec}{\mathbf{w}}

\newcommand{\Ocal}{\mathcal{O}}
\newcommand{\Om}{\Omega}
\newcommand{\ds}{\displaystyle}
\newcommand{\pa}{\partial}
\newcommand{\ul}{\underline}

\DeclareMathOperator{\rotv}{\mathbf{curl}}
\DeclareMathOperator{\grad}{\mathbf{grad}}
\DeclareMathOperator{\Grad}{\mathbf{Grad}}
\DeclareMathOperator{\Deltavec}{\mathbf{\Delta}}

\let\div=\divergence
\newtheorem{thm}{Theorem}

\newdefinition{rmk}{Remark}
\newtheorem{defi}{Definition}
\newtheorem{prop}{Proposition}
\newtheorem{cor}{Corollary}
\numberwithin{equation}{section}
\journal{Computers \& Mathematics with Applications}
\definecolor{cyan}{gray}{0}
\begin{document}
\begin{frontmatter}
\title{Explicit $T$-coercivity for the Stokes problem: a coercive finite element discretization\tnoteref{Title}}
\tnotetext[Title]{Explicit $T$-coercivity for Stokes}
\author[ENSTA]{Patrick Ciarlet~Jr}
\ead{patrick.ciarlet@ensta-paris.fr}
\affiliation[ENSTA]{organization={POEMS, CNRS, INRIA, ENSTA Paris, Institut Polytechnique de Paris},
            addressline={828 Boulevard
des Maréchaux}, 
            postcode={91120},
            city={Palaiseau},
            country={France}}
\author[CEA]{Erell Jamelot\corref{CorRef}}
\ead{erell.jamelot@cea.fr}
\cortext[CorRef]{Corresponding author: Erell Jamelot}
\affiliation[CEA]{organization={Universit\'e Paris-Saclay, CEA, Service de Thermo-hydraulique et de M\'ecanique des Fluides},
            postcode={91191},
            city={Gif-sur-Yvette cedex},
            country={France}}
\begin{abstract}
Using the $T$-coercivity theory as advocated in [Chesnel, Ciarlet, $T$-coercivity and continuous Galerkin methods: application to transmission problems with sign changing coefficients (2013)], we propose a new variational formulation of the Stokes problem which does not involve nonlocal operators. With this new formulation, unstable finite element pairs are stabilized. In addition, the numerical scheme is easy to implement, and a better approximation of the velocity and the pressure is observed numerically when the viscosity is small.
\end{abstract}
\end{frontmatter}
\input{CiJa2x-StokesGrad}
\section*{Acknowledgements}
E. Jamelot would like to thank CEA SIMU/SITHY project. 
\bibliographystyle{elsarticle-harv} 
\bibliography{CiJa2x-StokesGrad}
\appendix
\section{New variational formulation without orthogonality} \label{sec-New-VF-without-orthog}
Let us briefly consider what would happen if we were to consider a nonlocal right inverse of the divergence operator, that is different from the one proposed in Proposition \ref{prop:diviso}. For example, with values in the whole of $\bH^1_0(\Om)$ (not restricted to $\bV^\perp$). Such an operator $q\mapsto\bar\vvec_q$ is considered in Section 3.1 of \cite{BoFa13}, in the case of a 2D domain with a smooth boundary.
For this operator, let $\bar C_{\div}$ be a constant such that:
\begin{equation}\label{eq:BAItilde}
\forall q\in L^2_{zmv}(\Om),\,\exists \bar\vvec_q\in\bH^1_0(\Om)\,|\,\div\bar\vvec_q= q\mbox{ and }
\|\bar\vvec_q\|_{\bH^1_0(\Om)}\leq \bar C_{\div}\|q\|_{L^2(\Om)}.
\end{equation}
In general, $\bar\vvec_q$ does not belong to $\bV^\perp$, cf. \cite{BoFa13}. 

Let ${\lambda}>\frac{1}{4}(\bar C_{\div})^2$ (cf. remark \ref{rmk:ChoiceOfLambda}). $T$-coercivity can be obtained as before, with the operator 
\begin{equation}\label{eq:opTtilde}
\left\{
\begin{array}{rcl}
\bar T\colon\Xcal&\rightarrow\R\\
(\vvec,q)&\mapsto&({\lambda}\vvec- \nu^{-1}\bar\vvec_q,-{\lambda}q)
\end{array},
\right.
\end{equation}
Define $\bar{a}_{{\lambda}}\,\left((\uvec',p'),(\vvec,q)\right)=\textcolor{cyan}{a}\left((\uvec',p'),\bar T(\vvec,q)\right)$. We have:
\begin{equation}\label{eq:FVa1bar}
\begin{array}{rcl}
\ds \bar{a}_{{\lambda}}\,\left((\uvec',p'),(\vvec,q)\right)&=&
\ds \nu\,{\lambda}(\uvec',\vvec)_{\bH^1_0(\Om)}
- (\uvec',\bar\vvec_{q})_{\bH^1_0(\Om)} \\
\\
&&-{\lambda}(p',\div\vvec)_{L^2(\Om)}+\nu^{-1}\,(p',q)_{L^2(\Om)}\\
\\
&&+{\lambda}(q,\div\uvec')_{L^2(\Om)}.
\end{array}
\end{equation}
However, for $\uvec'\in\bV$ (that is for solutions $(\uvec',p)$ to Problem \eqref{eq:Stokes-FVa1tilde}), the term $(\uvec',\bar\vvec_q)_{\bH^1_0(\Om)}$ can no longer be removed from the expression (\ref{eq:FVa1tilde}) of the bilinear form $\bar{a}_{{\lambda}}$  \textcolor{cyan}{because orthogonality is lost}. As a matter of fact one has
\[ (\uvec',\bar\vvec_q)_{\bH^1_0(\Om)} = (\rotv\uvec',\rotv\bar\vvec_q)_{\bL^2(\Om)}.\]
Hence, in the non-orthogonal case, \textcolor{cyan}{one has to choose an ad hoc constant ${\lambda}>\frac{1}{4}(\bar C_{\div})^2$} to ensure $T$-coercivity (cf. Proposition~\ref{pro:a0Tcoer}). 
Similarly in the expression of the right-hand side $\bar\ell_{{\lambda}}$.
\begin{prop}\label{pro:fvq_nonorthog}
Let $\fvec'\in\bH^{-1}(\Om)$ be decomposed as in (\ref{eq:H-1}). Given, $q\in L^2_{zmv}(\Om)$, let $\bar\vvec_q\in\bH^1_0(\Om)$ be defined by (\ref{eq:BAItilde}). We have:
\begin{equation}\label{eq:fvq_nonorthog}
\langle\fvec',\bar\vvec_q\rangle_{\textcolor{cyan}{\bH^{-1}(\Om),\bH^1_0(\Om)}} 
 = -(z_{\fvec'},q)_{L^2(\Om)} \textcolor{cyan}{+ (\rotv\wvec_{\fvec'},\rotv\bar\vvec_q)_{\bL^2(\Om)}}.
\end{equation}
\end{prop}
To conclude \textcolor{cyan}{on the use of explicit $T$-coercivity in} the non-orthogonal case, we observe that the problem, if split as in (\ref{eq:StokesFV2}), leads to a more intricate variational formulation, which reads
\[ 
\left\{\begin{array}{l}
\mbox{Find $(\uvec,p)\in\Xcal$ s.t. fo r all $(\vvec,q)\in\Xcal$}\cr
(i)'\quad \nu \,{\lambda}(\uvec,\vvec)_{\bH^1_0(\Om)}- {\lambda}(p,\div\vvec)_{L^2(\Om)} = {\lambda}\langle\fvec,\vvec\rangle_{\textcolor{cyan}{\bH^{-1}(\Om),\bH^1_0(\Om)}},\\
(ii)'\ -(\uvec,\bar\vvec_{q})_{\bH^1_0(\Om)}+{\lambda}(q,\div\uvec)_{L^2(\Om)} + \nu^{-1}\,(p,q)_{L^2(\Om)}\\
\hskip 20truemm = \nu^{-1}\,(z_{\fvec},q)_{L^2(\Om)} \textcolor{cyan}{- \nu^{-1}\, (\rotv\wvec_{\fvec},\rotv\bar\vvec_q)_{\bL^2(\Om)}}.
\end{array}\right. \]
Regarding the right-hand side in $(ii')$, it requires the knowledge of both $z_{\fvec}$ and of $\wvec_{\fvec}$ (or of suitable approximations): a {\em two-step} procedure could be used as before (see the end of \S\ref{ss-sec-with_orthog}). One needs to evaluate also $(\uvec,\bar\vvec_{q})_{\bH^1_0(\Om)} = (\rotv\uvec,\rotv\bar\vvec_q)_{\bL^2(\Om)}$ in the left-hand side of $(ii')$. The latter part requires some knowledge of the nonlocal right inverse of the divergence operator, and it does not seem to be achievable numerically at a reasonable cost.

\section{\textcolor{cyan}{Nonhomogeneous Dirichlet boundary conditions}}\label{sec:CLDNH}
We now solve the classical incompressible Stokes model, with non zero Dirichlet boundary conditions. Hence, we focus on Problem \eqref{eq:Stokes}$_{NH}$ with $g=0$. The variational formulation 
reads:\\
Find $(\uvec,p)\in\bH^1(\Om)\times L^2_{zmv}(\Om)$ such that
\begin{equation}\label{eq:Stokes-NH}
\left\{\begin{array}{rcll}
\nu(\uvec,\vvec)_{\bH^1_0(\Om)}-(p,\div\vvec)_{L^2(\Om)}&=&\langle\fvec,\vvec\rangle_{\textcolor{cyan}{\bH^{-1}(\Om),\bH^1_0(\Om)}}&\forall\vvec\in\bH^1_0(\Om),\\
(q,\div\uvec)_{L^2(\Om)}&=&0&\forall q\in L^2_{zmv}(\Om),
\\
\uvec&=&\gvec\mbox{ on }\pa\Om,
\end{array}\right.
\end{equation}
where $\gvec\in\bH^\frac12(\pa\Om)$ is some boundary data such that $\gvec\cdot\nvec\in L^2_{zmv}(\pa\Om)$. Let $\uvec_{\gvec}\in\bH^1(\Om)$ be such that $\uvec_{\gvec}=\gvec\mbox{ on }\pa\Om$ and $-\Delta\uvec_{\gvec}=0$ in $\Om$. Then consider $(\uvec_0,p)\in\bH^1_0(\Om)\times L^2_{zmv}(\Om)$ the unique solution to \eqref{eq:Stokes}$_{H}$ with data $(\fvec+\nu\Delta\uvec_{\gvec},-\div\uvec_{\gvec})$:\\
Find $(\uvec_0,p)\in\Xcal$ such that for all $(\vvec,q)\in\Xcal$
\begin{equation}\label{eq:Stokes0-FV-NH}
\textcolor{cyan}{a}((\uvec_0,p),(\vvec,q))=
\langle\fvec,\vvec\rangle_{\textcolor{cyan}{\bH^{-1}(\Om),\bH^1_0(\Om)}}+(q,\div\uvec_{\gvec})_{L^2(\Om)},
\end{equation}
where we use that $(\uvec_{\gvec},\vvec)_{\bH^1_0(\Om)}=0$. Defining $\uvec=\uvec_{\gvec}+\uvec_0\in\bH^1(\Om)$, we find that $(\uvec,p)$ is solution to Problem \eqref{eq:Stokes-NH}. Uniqueness and continuous dependence with respect to the data are easily obtained. The next step is to replace $(\vvec,q)$ by $T((\vvec,q))=({\lambda}\vvec-\nu^{-1}\vvec_q,-{\lambda}q)$. One finds by orthogonality that for all $(\vvec,q)\in\Xcal$:
\[
\begin{array}{r}\ds
\nu{\lambda}(\uvec_0,\vvec)_{\bH^1_0(\Om)}-{\lambda}(p,\div\vvec)_{L^2(\Om)}+\nu^{-1}(p,q)_{L^2(\Om)}+{\lambda}(q,\div\uvec_0)_{L^2(\Om)}
\\=\ds
{\lambda}\langle\fvec,\vvec\rangle_{\textcolor{cyan}{\bH^{-1}(\Om),\bH^1_0(\Om)}}-\nu^{-1}\langle\fvec,\vvec_q\rangle_{\textcolor{cyan}{\bH^{-1}(\Om),\bH^1_0(\Om)}}
-{\lambda}(q,\div\uvec_{\gvec})_{L^2(\Om)}.
\end{array}
\]
Hence, a new variational formulation for nonhomogeneous boundary Dirichlet conditions is:\\
Find $(\uvec,p)\in\bH^1(\Om)\times L^2_{zmv}(\Om)$ such that for all $(\vvec,q)\in\Xcal$
\begin{equation}\label{eq:Stokes-FV-NH}
\begin{array}{c}
\nu{\lambda}(\uvec,\vvec)_{\bH^1_0(\Om)}-{\lambda}(p,\div\vvec)_{L^2(\Om)}+\nu^{-1}(p,q)_{L^2(\Om)}+{\lambda}(q,\div\uvec)_{L^2(\Om)}
\\=\ds
{\lambda}\langle\fvec,\vvec\rangle_{\textcolor{cyan}{\bH^{-1}(\Om),\bH^1_0(\Om)}}+\nu^{-1}(z_{\fvec},q)_{L^2(\Om)}.
\end{array}
\end{equation}
This is the same variational formulation as the one with homogeneous Dirichlet boundary conditions, except that $\uvec\in\bH^1(\Om)$ and $\uvec_{|\pa\Om}=\gvec$. \\

Regarding the numerical algorithms, if $z_{\fvec}$ is known, one solves a linear system like
\begin{equation}\label{eq:StokesFV2h-Exact_NH}
\left\{
\begin{array}{rcl}
\multicolumn{3}{l}{\mbox{Find }(\ul{U},\ul{P})\in\R^{N_u}\times\R^{N_p}\mbox{ such that:}}\\
\nu{\lambda}\bbA\,\ul{U}-{\lambda}\bbB^T\ul{P}&=&{\lambda}\ul{F}-\nu{\lambda}\bbA\,\ul{U}_{NH}\\
{\lambda}\bbB\,\ul{U}+\nu^{-1}\bbM\,\ul{P}&=&\nu^{-1}\bbM\,\ul{Z} -{\lambda}\bbB\,\ul{U}_{NH}
\end{array}
\right.
\end{equation}
where, in the right-hand side, $\ul{U}_{NH}\in\R^{N_u}$ accounts for $\uvec_{\gvec}$. This is completely similar to (\ref{eq:StokesFV2h-Exact}). While, if $z_{\fvec}$ is not known, starting from an initial guess $\ul{P}^{-1}\in\R^{N_p}$, for $n=0,1,\cdots$, one solves linear systems like
\begin{equation}\label{eq:StokesFV2h-Iterative_NH}
\left\{
\begin{array}{rcl}
\multicolumn{3}{l}{\mbox{Find }(\ul{U}^n,\ul{P}^n)\in\R^{N_u}\times\R^{N_p}\mbox{ such that:}}\\
\nu{\lambda}\bbA\,\ul{U}^n-{\lambda}\bbB^T\ul{P}^n&=&{\lambda}\ul{F}_u -\nu{\lambda}\bbA\,\ul{U}_{NH} \\
{\lambda}\bbB\,\ul{U}^n+\nu^{-1}\bbM\,\ul{P}^n&=&\nu^{-1}\bbM\,\ul{P}^{n-1} - {\lambda}\bbB\,\ul{U}_{NH}
\end{array}
\right..
\end{equation}
Interestingly, one recovers the same results as those of Theorem~\ref{thm-cv_algo}, because one finds identical iterating matrices, cf. (\ref{eq:deltaP}) and (\ref{eq:Yn}). \\

The (conforming) discretization of Problem \eqref{eq:Stokes-FV-NH} with $\bP^k\times P^{k-1}_{disc}$ finite element reads:
\begin{equation}
\label{eeq:Stokes-FV-NHh}
\left\{\begin{array}{l}
\mbox{Find }(\uvec_h,p_h)\in\bV_h^k\times Q_h^{k-1}\mbox{ s.t. for all }(\vvec_h,q_h)\in\bV_{0,h}^k\times Q_h^{k-1}
\\
\nu\,{\lambda}(\uvec_h,\vvec_h)_{\bH^1_0(\Om)}-{\lambda}(p_h,\div\vvec_h)_{L^2(\Om)} = {\lambda}\langle\fvec,\vvec_h\rangle_{\textcolor{cyan}{\bH^{-1}(\Om),\bH^1_0(\Om)}},\\
{\lambda}(q_h,\div\uvec_h)_{L^2(\Om)} + \nu^{-1}\,(p_h,q_h)_{L^2(\Om)} = \nu^{-1}\,(z_{\fvec},q_h)_{L^2(\Om)},\\
\uvec_{h|\pa\Om}=(\Pi_{h,cg}\uvec_{\gvec})_{|\pa\Om}.
\end{array}\right.
\end{equation}
For nonhomogeneous boundary Dirichlet conditions, one can study the error estimates by introducing $\uvec_{0,h}=\uvec_h-\Pi_{h,cg}\uvec_{\gvec}\in\bV_{0,h}^k$ and using the results that have been obtained for the homogeneous boundary conditions. First, in the estimate \eqref{eq:estim-p-ptilde} of Proposition~\ref{pro:zfh}, we change the term $\|z_{\fvec}-\tilde{z}_{\fvec,h}\|_{L^2(\Om)}$ in the right-hand side of the two equations into $\|z_{\fvec}-\tilde{z}_{\fvec,h}\|_{L^2(\Om)}+{\lambda}\nu\|\div(\uvec_{\gvec}-\Pi_{h,cg}\uvec_{\gvec})\|_{L^2(\Om)}$. Second, in the estimate \eqref{eq:estim-p-ptildeh} of Theorem \ref{th:zfh}, we add $\sqrt{{\lambda}}\|\div(\uvec_{\gvec}-\Pi_{h,cg}\uvec_{\gvec})\|_{L^2(\Om)}$ to the right-hand side.
\end{document}

%% file: CiJa2x-StokesGrad.tex
\section{Introduction}
The Stokes problem describes the steady state of incompressible Newtonian flows. They are derived from the Navier–Stokes equations \cite{GiRa86}. With regard to numerical analysis, the study of {the} Stokes problem helps to build an appropriate approximation of the Navier–Stokes equations. 
We propose here to write a new variational formulation of {the} Stokes problem using the $T$-coercivity theory, following Section 2 in \cite{ChCi13}. 

In Section \ref{sec:Tcoer}, we recall the $T$-coercivity theory as written in \cite{Ciar12,ChCi13}. 

\textcolor{cyan}{In Section \ref{sec:Stokes} we apply this theory to the generalized Stokes Problem, that is possibly with $\div\uvec\ne0$
In particular, solving this problem with homogeneous boundary condition helps solving {the} incompressible Stokes problem with inhomogeneous Dirichlet boundary conditions. We prove basic $T$-coercivity by finding} an operator $T$ such that the global variational formulation of the \textcolor{cyan}{generalized Stokes Problem} is $T$-coercive. 

\textcolor{cyan}{From that point on, we assume that $\div\uvec=0$ and solve the incompressible Stokes problem}. In Section \ref{sec:NVF}, we build and analyse a new variational formulation \textcolor{cyan}{with this operator $T$ (explicit $T$-coercivity). Then in Section~\ref{sec:algo} we propose a numerical algorithm based on the new variational formulation and, in the Section after, we introduce discretizations and we study more specifically stability properties. As a particular case, it is shown that the unstable finite element pair $\bP^1\times P^0$ yields stability}.
Finally, we provide some numerical experiments in Section~\ref{sec:resu} to illustrate our points. 

Some concluding remarks are made in Section~\ref{sec-conclusion}.
\section{$T$-coercivity}\label{sec:Tcoer}
We recall here the $T$-coercivity theory as written in \cite{Ciar12,ChCi13}. Consider the variational problem, where $V$ and $W$ are two Hilbert spaces and 
\textcolor{cyan}{$\ell\in W'$:
\begin{equation}
\label{eq:var-pb}
\mbox{Find }u\in V\mbox{ such that }\forall w\in W\mbox{, }a(u,w)= \ell(w).
\end{equation}
}
Classically, we know that Problem \eqref{eq:var-pb} is well-posed if $a(\cdot,\cdot)$ satisfies  the stability and the solvability conditions of the so-called Banach–Ne\v{c}as–Babu\v{s}ka (BNB) Theorem (see \textcolor{cyan}{e.g.} \cite[Theorem 25.9]{ErGu21-II}). For some models, one can also prove the well-posedness using the $T$-coercivity theory. Whereas the BNB theorem relies on an abstract inf–sup condition, $T$-coercivity uses explicit inf–sup operators. We refer for instance \textcolor{cyan}{to \cite{BoCZ10,NiVe11,BoCC12,ChCi13,BoCC14b,BoCC18} for problems with sign-changing coefficients}, \textcolor{cyan}{to \cite{BuCS02,BuCh03,Unge21} for integral equations, to \cite{BoCH11} for interior transmission problems}, to \cite{Hipt02,Buff05,Ciar12,HoNa15,Hall21} for Helmholtz-like problems, to \cite{JaCi13,CiJK17,Gire18} for the neutron diffusion equation, and to \cite{CiJa2x} for the magnetostatic problem. 
\begin{defi}[\textcolor{cyan}{basic $T$-coercivity}]\label{def:Tcoer}
Let $V$ and $W$ be two Hilbert spaces and $a(\cdot,\cdot)$ be a continuous bilinear form over $V\times W$. It is $T$-coercive if
\begin{equation}\label{eq:Tcoer}
\exists T\in\Lcal(V,W)\mbox{, bijective, }\exists\alpha>0\mbox{, }\forall v\in V\mbox{, }|a(v,Tv)|\geq\alpha\|v\|_V^2.
\end{equation}
\end{defi}
\textcolor{cyan}{Obviously, if (\ref{eq:Tcoer}) holds for some $T\in\Lcal(V,W)$, then $T$ is injective.}
\begin{thm}[well-posedness \cite{Ciar12,ChCi13}]\label{thm:WellPosed}
Let $a(\cdot,\cdot)$ be a continuous bilinear form on $V\times W$. \textcolor{cyan}{The Problem \eqref{eq:var-pb} is well-posed  if, and only if, the form $a(\cdot,\cdot)$ is $T$-coercive in the sense of definition~\ref{def:Tcoer}}.
\end{thm}
\textcolor{cyan}{Two extensions of $T$-coercivity are proposed in \cite[\S2.3.2]{ChCi13}. The first one is, when one is interested in solving the discretized problem, to mimic the design of the operator $T$ at the discrete level to obtain a uniform discrete inf-sup condition: we call it the {\em discrete $T$-coercivity}. Among the above cited references, this is developed in \cite{BuCS02,Hipt02,BoCZ10,Ciar12,ChCi13,HoNa15,CiJK17,BoCC18,Hall21,Unge21,Jame24}. The second one is the use of the operator $T$ directly in the variational formulation: an equivalent variational formulation is obtained, which reads
\begin{equation}
\label{eq:var-pb-Texplicit}
\mbox{Find }u\in V\mbox{ such that }\forall v\in V\mbox{, }a(u,Tv)= \ell (Tv).
\end{equation}
By design, $a(\cdot,T\cdot)$ is a continuous bilinear form, coercive on $V\times V$. We call it the {\em explicit $T$-coercivity}, see \cite{CiJa2x}.}
\section{\textcolor{cyan}{The} Stokes problem}\label{sec:Stokes}
Let $\Om$ be a connected bounded domain of $\R^d$, $d=2,\,3$, with a polygonal $(d=2)$ or Lipschitz polyhedral $(d=3)$ boundary $\pa\Om$. We consider \textcolor{cyan}{the generalized} Stokes problem:
\begin{equation}\label{eq:Stokes}
\mbox{Find }(\uvec,p)\mbox{ such that }
\begin{array}{rcl}
-\nu\Deltavec\uvec+\grad p &=&\fvec,\\
\div\uvec&=&\textcolor{cyan}{g},
\end{array}
\end{equation}
with Dirichlet boundary conditions for $\uvec$ and a normalization condition for $p$: $\ds\int_\Om p=0$. \textcolor{cyan}{If the Dirichlet boundary conditions are homogeneous, we write (\ref{eq:Stokes})$_{H}$, and (\ref{eq:Stokes})$_{NH}$ else.}

The vector field $\uvec$ represents the velocity of the fluid and the scalar field $p$ represents its pressure divided by the fluid density which is supposed to be constant. The first equation of \eqref{eq:Stokes} corresponds to the momentum balance equation and the second one corresponds to the conservation of the mass. The constant parameter $\nu>0$ is the kinematic viscosity of the fluid.\\ 
Let us provide some definition and reminders. Let us set $\bL^2(\Om):=(L^2(\Om))^d$, $\bH^1_0(\Om):=(H^1_0(\Om))^d$, $\bH^{-1}(\Omega):=(H^{-1}(\Om))^d$ \textcolor{cyan}{the dual space of $\bH^1_0(\Om)$} and $L^2_{zmv}(\Om):=\{q\in L^2(\Om)\,|\,\int_{\Om}q=0\}$, \textcolor{cyan}{$L^2_{zmv}(\pa\Om):=\{q\in L^2(\pa\Om)\,|\,\int_{\pa\Om}q=0\}$}.\\
\textcolor{cyan}{The natural function space for the velocity is $\bH^1(\Om)$ and, if homogeneous boundary conditions are prescribed, it is $\bH^1_0(\Om)$, while for the pressure the natural function space is $L^2_{zmv}(\Om)$. The data is $(\fvec,g)$.} The vector field $\fvec\in\bH^{-1}(\Om)$ represents a body force divided by the fluid density. \textcolor{cyan}{The scalar field $g\in L^2_{zvm}(\Om)$ is some abstract data that will be useful for the analysis of the case of non homogeneous Dirichlet boundary conditions.}

Let us first recall the Poincaré-Steklov inequality:
\begin{equation}\label{eq:Poincare}
\exists C_{PS}>0\,|\,\forall v\in H^1_0(\Om),\quad \|v\|_{L^2(\Om)}\leq C_{PS}\|\grad v\|_{\bL^2(\Om)}.
\end{equation}
Thanks to this result,  in $H^1_0(\Om)$ the semi-norm is equivalent to the full norm, so \textcolor{cyan}{we choose an inner} product  equal to $(v,w)_{H^1_0(\Om)}=(\grad v,\grad w)_{\bL^2(\Om)}$ \textcolor{cyan}{with associated norm} $\|v\|_{H^1_0(\Om)}=\|\grad v\|_{\bL^2(\Om)}$. Let $\vvec,\,\wvec\in\bH^1_0(\Om)$. We denote by $(v_i)_{i=1}^d$ (resp.\ $(w_i)_{i=1}^d$) the components of $\vvec$ (resp.\ $\wvec$), and we set $\Grad\vvec=(\pa_j v_i)_{i,j=1}^d\in\bbL^2(\Om)$, where $\bbL^2(\Om):=[L^2(\Om)]^{d\times d}$. We have: 
\[
(\Grad\vvec,\Grad\wvec)_{\bbL^2(\Om)}=(\vvec,\wvec)_{\bH^1_0(\Om)}=\ds\sum_{i=1}^d(v_i,w_i)_{H^1_0(\Om)}\]
and:
\[
\|\vvec\|_{\bH^1_0(\Om)}=\ds\left(\sum_{j=1}^d\|v_j\|_{H^1_0(\Om)}^2\right)^{1/2}=\|\Grad\vvec\|_{\bbL^2(\Om)}.
\]
Let us set $\bV:=\left\{\vvec\in\bH^1_0(\Om)\,|\,\div\vvec=0\right\}$. The vector space $\bV$ is a closed subset of $\bH^1_0(\Om)$. We denote by $\bV^\perp$ the orthogonal \textcolor{cyan}{complement} of $\bV$ in $\bH^1_0(\Om)$. We recall that \textcolor{cyan}{there exists a (nonlocal) right inverse of the divergence operator}:
\begin{prop}\label{prop:diviso}(\cite[Corollary I.2.4]{GiRa86})
The operator $\div\colon\bH^1_0(\Om)\rightarrow L^2(\Om)$ is an isomorphism of $\bV^{\perp}$ onto $L^2_{zmv}(\Om)$.
We call $C_{\div} \ge 1$ the smallest constant such that:
\begin{equation}\label{eq:BAI}
\forall p\in L^2_{zmv}(\Om),\,\exists! \tilde\vvec_p\in\bV^{\perp}\,|\,\div\tilde\vvec_p=p\mbox{ and }
\|\tilde\vvec_p\|_{\bH^1_0(\Om)}\leq C_{\div}\|p\|_{L^2(\Om)}.
\end{equation}
\end{prop}
In the \textcolor{cyan}{statement of Proposition~\ref{prop:diviso}}, we note that since
\[ \forall \vvec\in \bH^1_0(\Om), \quad \|\vvec\|_{\bH^1_0(\Om)}^2 = \|\rotv\vvec\|_{\bL^2(\Om)}^2 + \|\div\vvec\|_{L^2(\Om)}^2, \]
one has necessarily $C_{\div} \ge 1$. 

The variational formulation of Problem \eqref{eq:Stokes}$_{H}$ reads:\\
Find $(\uvec,p)\in\bH^1_0(\Om)\times L^2_{zmv}(\Om)$ such that
\begin{equation}\label{eq:Stokes-FV}
\begin{array}{rcll}
\nu(\uvec,\vvec)_{\bH^1_0(\Om)}-(p,\div\vvec)_{L^2(\Om)}&=&\langle\fvec,\vvec\rangle_{\textcolor{cyan}{\bH^{-1}(\Om),\bH^1_0(\Om)}}&\forall\vvec\in\bH^1_0(\Om)~;\\
(q,\div\uvec)_{L^2(\Om)}&=&\textcolor{cyan}{(g,q)_{L^2(\Om)}}&\forall q\in L^2_{zmv}(\Om).
\end{array}
\end{equation}
\textcolor{cyan}{This is a saddle point-problem. Using classical theory, one proves} that Problem \eqref{eq:Stokes-FV} is well-posed with the help of Poincaré-Steklov inequality \eqref{eq:Poincare} and Proposition \ref{prop:diviso}. Check for instance the proof of \cite[Theorem I.5.1]{GiRa86}.

Let us set $\Xcal:=\bH^1_0(\Om)\times L^2_{zmv}(\Om)$ which is a Hilbert space which we endow with the following norm:
\begin{equation}\label{eq:normX}
 \|(\vvec,q)\|_{\Xcal,\nu}=\left(\|\vvec\|_{\bH^1_0(\Om)}^2+\nu^{-2}\,\|q\|_{L^2(\Om)}^2\right)^{1/2}.
\end{equation}
\textcolor{cyan}{This norm is chosen to account for physical phenomena due to small viscosity. Typically, if the result of a physics experiment is $(\uvec,p)$, the ratio of the two components of the norm $\|(\uvec,p)\|_{\Xcal,\nu}$, respectively equal to $\|\uvec\|_{\bH^1_0(\Om)}$ and $\nu^{-1}\|p\|_{L^2(\Om)}$, varies linearly with $\nu$.} \\
We define the following continuous symmetric bilinear form:
\begin{equation}\label{eq:aform}
\left\{\begin{array}{rcl}\textcolor{cyan}{a}\colon
\Xcal\times\Xcal&\rightarrow&\R\\
(\uvec',p')\times(\vvec,q)&\mapsto&\nu(\uvec',\vvec)_{\bH^1_0(\Om)}-(p',\div\vvec)_{L^2(\Om)}-(q,\div\uvec')_{L^2(\Om)}
\end{array}.\right.
\end{equation}
We also define the linear continuous form:
\begin{equation}\label{eq:ell0}
\left\{\begin{array}{rcl}
\textcolor{cyan}{\ell}\colon\Xcal&\rightarrow&\R\\
(\vvec,q)&\mapsto&\langle\fvec,\vvec\rangle_{\textcolor{cyan}{\bH^{-1}(\Om),\bH^1_0(\Om)}}\textcolor{cyan}{-(g,q)_{L^2(\Om)}}
\end{array}.\right.
\end{equation}
We can write Problem \eqref{eq:Stokes}$_{H}$ in an equivalent way as follows:
\begin{equation}\label{eq:Stokes-FVa}
\mbox{Find }(\uvec,p)\in\Xcal\mbox{ s.t.}\quad \textcolor{cyan}{a}\left((\uvec,p),(\vvec,q)\right)=\textcolor{cyan}{\ell}((\vvec,q))\quad\forall(\vvec,q)\in\Xcal.
\end{equation}

Let us prove that Problem \eqref{eq:Stokes-FVa} is well-posed using basic $T$-coercivity.
\begin{prop}\label{pro:a0Tcoer}
The bilinear form $\textcolor{cyan}{a}(\cdot,\cdot)$ is $T$-coercive.
\end{prop}
\begin{proof}
We follow here the proof given in \cite{BaCi22,BoCi22}. Let us consider $(\uvec',p')\in\Xcal$ and let us build $(\vvec^\star,q^\star)=T(\uvec',p')\in\Xcal$ satisfying \eqref{eq:Tcoer} (with $\textcolor{cyan}{V=W=\Xcal}$). We need three main steps.
\begin{itemize}
\item[1.] According to Proposition \ref{prop:diviso}, there exists $\tilde{\vvec}_{p'}\in\bV^\perp$ such that:
\begin{equation}
\label{eq:diviso}
\div\tilde{\vvec}_{p'}=p'\mbox{ in }\Om\mbox{ and }\|\tilde{\vvec}_{p'}\|_{\bH^1_0(\Om)}\leq C_{\div}\,\|p'\|_{L^2(\Om)}.
\end{equation}
Let us set $(\vvec^\star,q^\star):=(\lambda\uvec'-\nu^{-1}\tilde{\vvec}_{p'},-\lambda p')$, with \textcolor{cyan}{$\lambda>0$ to be fixed}. We obtain:
\begin{equation}\label{eq:coer1}
\textcolor{cyan}{a}\left((\uvec',p'),(\vvec^\star,q^\star)\right)=
\nu\,\lambda\|\uvec'\|_{\bH^1_0(\Om)}^2+\nu^{-1}\,\|p'\|_{L^2(\Om)}^2-\,(\uvec',\tilde{\vvec}_{p'})_{\bH^1_0(\Omega)}.
\end{equation}
\item[2.] In order to bound the last term of \eqref{eq:coer1}, we use the Young inequality and then inequality \eqref{eq:diviso}, and it follows that for all $\eta>0$:
\begin{equation}\label{eq:GvGvp}
(\uvec',\tilde{\vvec}_{p'})_{\bH^1_0(\Omega)}\leq
\ds\frac{\eta}{2}\|\uvec'\|_{\bH^1_0(\Omega)}^2+\frac{\eta^{-1}}{2}\,(C_{\div})^2\,\|p'\|_{L^2(\Om)}^2.
\end{equation}
\item[3.] 
Using the bound \eqref{eq:GvGvp} in \eqref{eq:coer1}:
\[
\textcolor{cyan}{a}\left((\uvec',p'),(\vvec^\star,q^\star)\right)\geq \left(\nu\,\lambda - \frac{\eta}{2} \right)\|\uvec'\|_{\bH^1_0(\Omega)}^2+ \left(\nu^{-1} - \frac{\eta^{-1}}{2}\,(C_{\div})^2\right)\|p'\|_{L^2(\Om)}^2.
\]
\textcolor{cyan}{We look for $\eta>0$ such that $2 \nu \lambda > \eta$ and $\eta > \frac{\nu}{2}\,(C_{\div})^2$, which amounts to requiring
\[  \lambda > \frac14\,(C_{\div})^2. \]}
\end{itemize}
According to the above, provided that $\lambda > \frac14\,(C_{\div})^2$, we have proved that the operator $T_{\lambda}\in\Lcal(\Xcal)$ defined by $T_{\lambda}\left((\uvec',p')\right)=(\lambda\uvec'-\nu^{-1}\tilde{\vvec}_{p'},-\lambda p')$ is such that:
\[ \exists \alpha_{\lambda}>0,\ \forall (\uvec',p')\in \Xcal, \textcolor{cyan}{a}\left((\uvec',p'),T_{\lambda}\left((\uvec',p')\right)\right)\geq \alpha_{\lambda}\,\|(\uvec',p')\|_{\Xcal,\nu}^2. \]
\textcolor{cyan}{As noted in Section~\ref{sec:Tcoer}}, the injectivity of the operator $T_{\lambda}$ \textcolor{cyan}{follows}. Given $(\vvec^\star,q^\star)\in\Xcal$, choosing $(\uvec',p')=((\lambda^{-1}\vvec^\star-\nu^{-1}\lambda^{-2}\tilde{\vvec}_{q^\star},-\lambda^{-1}q^\star)$ yields $T_{\lambda}((\uvec',p'))=(\vvec^\star,q^\star)$. Hence, the operator $T_{\lambda}\in\Lcal(\Xcal)$ is bijective.
\end{proof}
\textcolor{cyan}{
\begin{rmk}\label{rmk:ChoiceOfLambda}
In the above proof, we established that the the bijective operator $T_{\lambda}\left((\uvec',p')\right)=(\lambda\uvec'-\nu^{-1}\tilde{\vvec}_{p'},-\lambda p')$ leads to $T$-coercivity as soon as $\lambda>\frac14(C_{\div})^2$. Observe that one can have even more flexibility in the choice of $T$ by choosing a different factor in front of $\tilde{\vvec}_{p'}$, and then choosing $\lambda$ accordingly.
\end{rmk}
}
We can now prove \textcolor{cyan}{the following result for the generalized Stokes problem:}
\begin{thm}
\label{thm:StokesWellPosed}
Problem \eqref{eq:Stokes-FVa} is well-posed, so it admits one and only one solution for any $(\fvec,g)\in\bH^{-1}(\Om)\times L^2_{zmv}(\Om)$. \textcolor{cyan}{Writing $\uvec=\uvec_0+\uvec_\perp$ with $\uvec_0\in\bV$ and $\uvec_\perp\in\bV^\perp$}, the solution is such that:
\begin{equation}\label{eq:WellPosed}
\forall\fvec\in\bH^{-1}(\Omega),\,\textcolor{cyan}{\forall g\in L^2_{zmv}(\Om)}\quad\left\{
\begin{array}{l}
\textcolor{cyan}{\|\uvec_\perp\|_{\bH^1_0(\Om)}\leq C_{\div}\|g\|_{L^2(\Om)}},\\
\|\textcolor{cyan}{\uvec_0}\|_{\bH^1_0(\Om)}\leq\nu^{-1}\,\|\fvec\|_{\bH^{-1}(\Omega)},\\
\|p\|_{L^2(\Om)}\leq C_{\div}\,\|\fvec\|_{\bH^{-1}(\Om)}\textcolor{cyan}{+\nu\,C_{\div}^2\,\|g\|_{L^2(\Om)}}.\\
\end{array}\right.
\end{equation}
\end{thm}
\begin{proof}
According to Proposition \ref{pro:a0Tcoer}, the continuous bilinear form $\textcolor{cyan}{a}(\cdot,\cdot)$ is $T$-coercive. Hence, according to Theorem \ref{thm:WellPosed}, Problem \eqref{eq:Stokes-FVa} is well-posed. Let us now derive \eqref{eq:WellPosed}. Consider $(\uvec,p)$ the unique solution of Problem \eqref{eq:Stokes-FVa}, \textcolor{cyan}{where $\uvec=\uvec_0+\uvec_\perp$ with $\uvec_0\in\bV$ and $\uvec_\perp\in\bV^\perp$}. Choosing $\vvec=0$, we obtain that \textcolor{cyan}{$\div\uvec_\perp=g$, so that $\|\uvec_\perp\|_{\bH^1_0(\Om)}\leq C_{\div}\|g\|_{L^2(\Om)}$}. Now, choosing $\vvec=\textcolor{cyan}{\uvec_0}$ and using \textcolor{cyan}{orthogonality}, we have: $\nu\,\|\textcolor{cyan}{\uvec_0}\|_{\bH^1_0(\Om)}^2=\langle\fvec,\textcolor{cyan}{\uvec_0}\rangle_{\textcolor{cyan}{\bH^{-1}(\Om),\bH^1_0(\Om)}}\leq\|\fvec\|_{\bH^{-1}(\Om)}\,\|\textcolor{cyan}{\uvec_0}\|_{\bH^1_0(\Om)}$, so that: $\|\textcolor{cyan}{\uvec_0}\|_{\bH^1_0(\Om)}\leq\nu^{-1}\,\|\fvec\|_{\bH^{-1}(\Om)}$. Next, we choose, in \eqref{eq:Stokes-FVa}, $q=0$ and $\vvec=-\tilde{\vvec}_{p}\in\bV^{\perp}$, where $\div\tilde{\vvec}_p= p$ and $\|\tilde{\vvec}_{p}\|_{\bH^1_0(\Om)}\leq C_{\div}\,\|p\|_{L^2(\Om)}$ (see Proposition \ref{prop:diviso}). Since $\textcolor{cyan}{\uvec_0}\in\bV$, 
it holds that $(\textcolor{cyan}{\uvec_0},\tilde{\vvec}_{p})_{\bH^1_0(\Om)}=0$. This gives:
\[
\begin{array}{rcl}
\|p\|_{L^2(\Om)}^2&=&(p,\div\tilde{\vvec}_p)_{L^2(\Om)} = -\langle\fvec, \tilde{\vvec}_{p}\rangle_{\textcolor{cyan}{\bH^{-1}(\Om),\bH^1_0(\Om)}}\textcolor{cyan}{+\nu\langle\Delta\uvec_\perp,\tilde{\vvec}_{p}\rangle_{\textcolor{cyan}{\bH^{-1}(\Om),\bH^1_0(\Om)}}},\\
&\leq&
\left(\|\fvec\|_{\bH^{-1}(\Om)}+\textcolor{cyan}{\nu\|\uvec_{\perp}\|_{\bH^1_0(\Om)}}\right)\,\|\tilde{\vvec}_{p}\|_{\bH^1_0(\Om)}\\
&\leq& C_{\div}\,\left(\|\fvec\|_{\bH^{-1}(\Om)}\textcolor{cyan}{+\nu\,C_{\div}\,\|g\|_{L^2(\Om)}}\right)\,\|p\|_{L^2(\Om)},
\end{array}
\]
so that: $\|p\|_{L^2(\Om)}\leq C_{\div}\,\|\fvec\|_{\bH^{-1}(\Om)}\textcolor{cyan}{+\nu\,C_{\div}^2\,\|g\|_{L^2(\Om)}}$.
\end{proof}
\textcolor{cyan}{
In the following sections, we will consider that $g=0$, i.e. we focus on the classical incompressible Stokes model.
}
\section{New variational {formulations}}\label{sec:NVF}
\textcolor{cyan}{
We solve the classical incompressible Stokes model, with $g=0$.}
\subsection{\textcolor{cyan}{Explicit} $T$-coercivity}
{Let $\lambda>\frac14(C_{\div})^2$. In remark \ref{rmk:ChoiceOfLambda}, we introduced the operator $T_{\lambda}\in\Lcal(\Xcal)$ defined by}:
\begin{equation}\label{eq:opT}
T_{\lambda}\colon\left\{
\begin{array}{rcl}
\Xcal&\rightarrow&\R\\
(\vvec,q)&\mapsto&(\lambda\vvec-\nu^{-1}\tilde{\vvec}_{q},-\lambda q)
\end{array},
\right.
\end{equation}
where $\tilde{\vvec}_{q}\in\bV^{\perp}$ is given by (\ref{eq:BAI}): $\div\tilde{\vvec}_{q}= q$ and $\|\tilde{\vvec}_q\|_{\bH^1_0(\Om)}\leq C_{\div}\,\|q\|_{L^2\Om)}$. 
We now write the variational formulation \eqref{eq:Stokes-FVa} with test function $T_{\lambda}\left((\vvec,q)\right)$ instead of $(\vvec,q)$. Let us define $\tilde{a}_{\lambda}\,\left((\uvec',p'),(\vvec,q)\right)=\textcolor{cyan}{a}\left((\uvec',p'),T(\vvec,q)\right)$. We have:
\begin{equation}\label{eq:FVa1tilde}
\begin{array}{rcl}
\ds \tilde{a}_{\lambda}\,\left((\uvec',p'),(\vvec,q)\right)&=&
\ds \nu\,\lambda(\uvec',\vvec)_{\bH^1_0(\Om)}
-\,(\uvec',\tilde{\vvec}_{q})_{\bH^1_0(\Om)}\\
\\
&&-\lambda(p',\div\vvec)_{L^2(\Om)}+\nu^{-1}\,(p',q)_{L^2(\Om)}\\
\\
&&+\lambda(q,\div\uvec')_{L^2(\Om)}.
\end{array}
\end{equation}
\textcolor{cyan}{The term with $\tilde{\vvec}_{q}$ requires some knowledge of the nonlocal right inverse of the divergence operator. We will see in the next section that it can be removed.} \\
According to Proposition \ref{pro:a0Tcoer}, we have the...
\begin{prop}
The bilinear form \eqref{eq:FVa1tilde} is coercive.
\end{prop}
Introducing $\ell_{\lambda}((\vvec,q)):=\textcolor{cyan}{\ell}(T_{\lambda}(\vvec,q))$, we can propose a first new variational formulation to Problem \eqref{eq:Stokes}$_{H}$, which reads:
\begin{equation}\label{eq:Stokes-FVa1tilde}
\mbox{Find }(\uvec,p)\in\Xcal\mbox{ s.t.}\quad \tilde{a}_{\lambda}\,\left((\uvec,p),(\vvec,q)\right)=\ell_{\lambda}((\vvec,q))\quad\forall(\vvec,q)\in\Xcal.
\end{equation}
Using the bijectivity of $T_{\lambda}$, it is obvious that (\ref{eq:Stokes-FVa1tilde}) is equivalent to \eqref{eq:Stokes-FVa}, so well-posedness of (\ref{eq:Stokes-FVa1tilde}) is a direct consequence of the well-posedness of \eqref{eq:Stokes-FVa}, and vice versa.

Let us derive an "explicit" expression of $\ell_{\lambda}((\vvec,q))$, which will prove useful later on. We recall that we have in the sense of distributions:
\begin{equation}
\label{eq:deltarotgrad}
-\Deltavec(\cdot)=\textcolor{cyan}{\rotv\rotv(\cdot)-\grad\div(\cdot)},
\end{equation}
\textcolor{cyan}{so that}
\begin{equation}\label{eq:H-1}
\begin{array}{c}
\forall\fvec'\in\bH^{-1}(\Om),\quad\exists! (z_{\fvec'},\wvec_{\fvec'})\in L^2_{zmv}(\Om)\times\bV\,|\,\\
\\
\fvec'\textcolor{cyan}{=\grad z_{\fvec'}-\Deltavec\wvec_{\fvec'}}=\grad z_{\fvec'}+\rotv\rotv\wvec_{\fvec'},
\end{array}
\end{equation}
where $(\wvec_{\fvec'},z_{\fvec'})$ solves \eqref{eq:Stokes}${}_H$ with $\nu=1$ and data \textcolor{cyan}{$\fvec=\fvec'$ and $g=0$}.

\begin{prop}\label{pro:fvq}
{Let $\fvec'\in\bH^{-1}(\Om)$. Given, $q\in L^2_{zmv}$, let $\tilde\vvec_q\in\bV^{\perp}$ be defined by (\ref{eq:BAI})}.
We have:
\begin{equation}\label{eq:fvq}
\langle\fvec',\tilde\vvec_q\rangle_{\textcolor{cyan}{\bH^{-1}(\Om),\bH^1_0(\Om)}} 
 = -(z_{\fvec'},q)_{L^2(\Om)}.
\end{equation}
\end{prop}
\begin{proof} Let $\fvec'$ be decomposed as in (\ref{eq:H-1}). 

{On the one hand,} integrating by parts twice and using \eqref{eq:deltarotgrad}, we get:
\begin{equation*}\label{eq:rotw}
\begin{array}{rcl}
-\langle\rotv\rotv\wvec_{\fvec'}
,\tilde\vvec_q\rangle_{\textcolor{cyan}{\bH^{-1}(\Om),\bH^1_0(\Om)}}&=&-\langle\rotv\rotv\tilde\vvec_q,\wvec_{\fvec'}\rangle_{\textcolor{cyan}{\bH^{-1}(\Om),\bH^1_0(\Om)}},\\
&=& \textcolor{cyan}{\langle\Deltavec\tilde\vvec_q,\wvec_{\fvec'}\rangle_{\textcolor{cyan}{\bH^{-1}(\Om),\bH^1_0(\Om)}}}- \langle\grad q,\wvec_{\fvec'}\rangle_{\textcolor{cyan}{\bH^{-1}(\Om),\bH^1_0(\Om)}},\\
&=& \textcolor{cyan}{-(\tilde\vvec_q,\wvec_{\fvec'})_{\bH^1_0(\Om)}} +  (q,\div\wvec_{\fvec'})_{L^2(\Om)},\\
&=&0,
\end{array}
\end{equation*}
resp.\ since {$\div\tilde\vvec_q=q$}, $\tilde\vvec_q\in\bV^{\perp}$, $\wvec_{\fvec'}\in\bV$, and $\div\wvec_{\fvec'}=0$. 

On the other hand, integrating by parts once, we have:
\begin{equation*}\label{eq:gradz}
\langle\grad z_{\fvec'},\tilde\vvec_q\rangle_{\textcolor{cyan}{\bH^{-1}(\Om),\bH^1_0(\Om)}}=-\,(z_{\fvec'},q)_{L^2(\Om)},
\end{equation*}
so the claim follows.
\end{proof}
\textcolor{cyan}{So, we have the}
\begin{cor}\label{cor_l_lambda} The right-hand-side $\ell_{\lambda}$ is equal to \[(\vvec,q)\mapsto \lambda\langle\fvec,\vvec\rangle_{\textcolor{cyan}{\bH^{-1}(\Om),\bH^1_0(\Om)}}+\nu^{-1}(z_{\fvec},q)_{L^2(\Om)}.\]
\end{cor}
 \begin{proof} 
\textcolor{cyan}{This is an obvious consequence of \eqref{eq:fvq}}.
\end{proof}
\subsection{New variational formulation using orthogonality}\label{ss-sec-with_orthog}
Going back to the original Problem \eqref{eq:Stokes}$_{H}$, we note that all solutions $(\uvec',p')$ to Problem \eqref{eq:Stokes-FVa1tilde} are such that $\uvec'\in\bV$. So, in the statement of Problem~(\ref{eq:Stokes-FVa1tilde}), the term $(\uvec',\tilde\vvec_q)_{\bH^1_0(\Om)}$ can be removed \textcolor{cyan}{by orthogonality} from the expression (\ref{eq:FVa1tilde}) of the bilinear form $\tilde{a}_{\lambda}$: as a matter of fact, according to (\ref{eq:BAI}), for all $q\in L^2_{zmv}(\Om)$, $(\uvec',\tilde\vvec_q)_{\bH^1_0(\Om)}=0$ because $\tilde\vvec_q\in\bV^\perp$. Interestingly, in this manner one improves the stability in the proof of Proposition~\ref{pro:a0Tcoer}, because the cross term vanishes in the expression (\ref{eq:coer1}): no treatment is required and the resulting bilinear form is coercive for all $\lambda>0$. So, \textcolor{cyan}{we propose to remove the second term in the expression (\ref{eq:FVa1tilde}) of the bilinear form with operator $T_{\lambda}$, that is, we} introduce the bilinear form $a_{\lambda}(\cdot,\cdot)$ such that:
\begin{equation}\label{eq:FVa1}
\left\{
\begin{array}{rcl}
a_{\lambda}\colon\quad\Xcal\times\Xcal&\rightarrow&\R\\
\left((\uvec',p'),(\vvec,q)\right)&\mapsto&\nu\,\lambda(\uvec',\vvec)_{\bH^1_0(\Om)}+\nu^{-1}\,(p',q)_{L^2(\Om)}
\\
\\
&&+\,\lambda\left[\,(q,\div\uvec')_{L^2(\Om)}-(p',\div\vvec)_{L^2(\Om)}\,\right]
\end{array}.
\right.
\end{equation}
\begin{prop}
The bilinear form \eqref{eq:FVa1} is coercive.
\end{prop}
\begin{proof}
We have:
\[\begin{array}{rcl}a_{\lambda}\left((\uvec',p'),(\uvec',p')\right)&=&\nu\,\lambda\|\uvec'\|_{\bH^1_0(\Om)}^2+\nu^{-1}\,\|p'\|^2_{L^2(\Om)},\\
&\geq&\ds\nu\,\min(1,\lambda)\,\|(\uvec',p')\|_{\Xcal,\nu}^2.
\end{array}
\]
\end{proof}
\textcolor{cyan}{With the help of explicit $T$-coercivity and using orthogonality}, we can now propose a second new variational formulation to Problem \eqref{eq:Stokes}$_{H}$, which reads:
\begin{equation}\label{eq:Stokes-FVa1}
\left\{
\begin{array}{l}
\mbox{Find }(\uvec,p)\in\Xcal\mbox{ such that }\cr
 a_{\lambda}\left((\uvec,p),(\vvec,q)\right)={\lambda}{\langle\fvec,\vvec\rangle_{\textcolor{cyan}{\bH^{-1}(\Om),\bH^1_0(\Om)}}+\nu^{-1}\,(z_{\fvec},q)_{L^2(\Om)}}\quad\forall(\vvec,q)\in\Xcal.
\end{array}
\right.
\end{equation}
\begin{thm}\label{th:equiv1}
\textcolor{cyan}{For all $\lambda>0$, Problem} \eqref{eq:Stokes-FVa1} is well-posed and is equivalent to Problem \eqref{eq:Stokes}$_{H}$ \textcolor{cyan}{with $g=0$}.
\end{thm}
\textcolor{cyan}{
By contrast with remark \ref{rmk:ChoiceOfLambda}, the result here holds for all $\lambda>0$.
}
\begin{proof} The bilinear form $a_{\lambda}(\cdot,\cdot)$ is continuous and coercive. Let $\fvec\in\bH^{-1}(\Om)$, and let $z_{\fvec}$ be given by (\ref{eq:H-1}), the linear form $\ell_{\lambda}(\,\cdot\,$) is continuous. According to Lax-Milgram Theorem, Problem \eqref{eq:Stokes-FVa1} is well-posed. It exists a unique solution $(\uvec,p)$ which depends continuously on the data. 

Regarding the equivalence with Problem \eqref{eq:Stokes}$_{H}$ \textcolor{cyan}{with $g=0$}, we already observed that solving (\ref{eq:Stokes-FVa1tilde}) is equivalent to solving \eqref{eq:Stokes-FVa}, and that both of them are equivalent to solving Problem \eqref{eq:Stokes}$_{H}$ \textcolor{cyan}{with $g=0$}. Then, if $(\uvec,p)$ solves Problem \eqref{eq:Stokes}$_{H}$ with data $\fvec\in\bH^{-1}(\Om)$ \textcolor{cyan}{and $g=0$}, one has in particular that $\uvec\in\bV$. Hence it follows from the above that $(\uvec,p)$ solves (\ref{eq:Stokes-FVa1}) with data $\fvec$ \textcolor{cyan}{and $g=0$}, resp.\ $z_{\fvec}$ given by (\ref{eq:H-1}). 

Conversely, assume that $(\uvec,p)$ solves (\ref{eq:Stokes-FVa1}), with data $\fvec\in\bH^{-1}(\Om)$, and $z_{\fvec}$ given by (\ref{eq:H-1}). Denote by $(\uvec^\dagger,p^\dagger)$ the solution to Problem \eqref{eq:Stokes}$_{H}$ with data $\fvec$. According to the above and by linearity, $(\uvec-\uvec^\dagger,p-p^\dagger)$ solves (\ref{eq:Stokes-FVa1}) with vanishing data, hence to is equal to $({\bf0},0)$ by uniqueness: in other words, $(\uvec,p)$ solves Problem \eqref{eq:Stokes}$_{H}$ with data $\fvec$.
\end{proof}
Notice that we can write problem \eqref{eq:Stokes-FVa1} with two equations as:
\begin{equation}
\label{eq:StokesFV2}
\left\{\begin{array}{l}
\mbox{Find $(\uvec,p)\in\Xcal$ s.t. for all $(\vvec,q)\in\Xcal$}\cr
(i)\quad \nu\,\lambda(\uvec,\vvec)_{\bH^1_0(\Om)}-\lambda(p,\div\vvec)_{L^2(\Om)} = \lambda\langle\fvec,\vvec\rangle_{\textcolor{cyan}{\bH^{-1}(\Om),\bH^1_0(\Om)}},\\
(ii)\ \lambda(q,\div\uvec)_{L^2(\Om)} + \nu^{-1}\,(p,q)_{L^2(\Om)} = \nu^{-1}\,(z_{\fvec},q)_{L^2(\Om)}.
\end{array}\right.
\end{equation}
This new variational formulation appears as a stabilized variational formulation, in the sense of {\S}II.1.2 in \cite{GiRa86}, pages 120-123. It can be solved once $z_{\fvec}$, or a suitable approximation of $z_{\fvec}$, is available. This suggests to use (\ref{eq:StokesFV2}) as a {\em post processing step} as follows:
\begin{itemize}
\item Compute $z_{\fvec}$ by solving numerically the Stokes problem (\ref{eq:Stokes-FV}) with $\nu=1$, and data $\fvec$\,;
\item Solve (\ref{eq:StokesFV2}) with the data $\fvec$ and the computed $z_{\fvec}$.
\end{itemize}
Notice that we can recover that the solution $\uvec$ to (\ref{eq:StokesFV2}) belongs to $\bV$ in a simple manner: let us split $\uvec=\uvec_0+\uvec_{\perp}$, where $(\uvec_0,\uvec_{\perp})\in\bV\times\bV^\perp$, so that $\div\uvec=\div\uvec_{\perp}$. Choosing $\vvec=\uvec_{\perp}$ in \eqref{eq:StokesFV2}-$(i)$, and $q=\div\uvec$ in \eqref{eq:StokesFV2}-$(ii)$:
\begin{equation}
\label{eq:StokesFV2-div0}
\left\{\begin{array}{ll}
(i)\quad\nu\,\|\uvec_{\perp}\|_{\bH^1_0(\Om)}^2-(p,\div\uvec)_{L^2(\Om)}=\langle\fvec,\uvec_{\perp}\rangle_{\textcolor{cyan}{\bH^{-1}(\Om),\bH^1_0(\Om)}}=-(z_{\fvec},\div\uvec)_{L^2(\Om)},\\
(ii)\quad\lambda\|\div\uvec\|^2_{L^2(\Om)}+\nu^{-1}(p,\div\uvec)_{L^2(\Om)}=\nu^{-1}(z_{\fvec},\div\uvec)_{L^2(\Om)}.
\end{array}\right.
\end{equation}
Summing \eqref{eq:StokesFV2-div0}-$(i)$ and \eqref{eq:StokesFV2-div0}-$(ii)$ times $\nu$, we obtain:
\[
\nu\,\left(\|\uvec_{\perp}\|_{\bH^1_0(\Om)}^2+\lambda\|\div\uvec\|^2_{L^2(\Om)}\right)=0.
\]
Hence, $\uvec_{\perp}=0$ and $\div\uvec=0$. \\ 

\textcolor{cyan}{In \ref{sec-New-VF-without-orthog}, we present another variational formulation \textcolor{cyan}{(not relying on orthogonality)} that can be obtained with the help of explicit $T$-coercivity. In \ref{sec:CLDNH}, we address the case of the classical incompressible Stokes model, with non zero Dirichlet boundary conditions.}
\section{\textcolor{cyan}{Numerical algorithms}}\label{sec:algo}
We solve the classical incompressible Stokes model, with $g=0$. We propose below two strategies depending on whether $z_{\fvec}$ is known or not. {Below, we consider homogeneous Dirichlet boundary conditions (see \ref{sec:CLDNH} for nonhomogeneous Dirichlet boundary conditions)}. \\
If $z_{\fvec}$ is known, one solves (\ref{eq:StokesFV2}) directly. If this is not the case, we propose to start from some approximate value of $z_{\fvec}$, and then to iterate by solving a series of problems like (\ref{eq:StokesFV2}). \\
So, after discretization, if $z_{\fvec}$ is known, one solves a linear system like
\begin{equation}\label{eq:StokesFV2h-Exact}
\left\{
\begin{array}{rcl}
\multicolumn{3}{l}{\mbox{Find }(\ul{U},\ul{P})\in\R^{N_u}\times\R^{N_p}\mbox{ such that:}}\\
\nu\lambda\bbA\,\ul{U}-\lambda\bbB^T\ul{P}&=&\lambda\ul{F}\\
\lambda\bbB\,\ul{U}+\nu^{-1}\bbM\,\ul{P}&=&\nu^{-1}\bbM\,\ul{Z}
\end{array}
\right.
\end{equation}
where $\ul{U}\in\R^{N_u}$ represents the discrete velocity, resp. $\ul{P}\in\R^{N_p}$ the discrete pressure, while $\ul{F}\in\R^{N_u}$ stands for $\fvec$, and $\ul{Z}\in\R^{N_p}$ stands for $z_{\fvec}$. Classically, $\ul{U}\mapsto (\bbA\ul{U}|\ul{U})$ measures discrete velocities, $\ul{P}\mapsto (\bbM\ul{P}|\ul{P})$ measures discrete pressures and, as a consequence, $\ul{U}\mapsto (\bbB\ul{U}|\bbM^{-1}\bbB\ul{U})$ measures the divergence of the discrete velocities. \\
On the other hand, if $z_{\fvec}$ is not known, starting from an initial guess $\ul{P}^{-1}\in\R^{N_p}$,
for $n=0,1,\cdots$, one solves linear systems like
\begin{equation}\label{eq:StokesFV2h-Iterative}
\left\{
\begin{array}{rcl}
\multicolumn{3}{l}{\mbox{Find }(\ul{U}^n,\ul{P}^n)\in\R^{N_u}\times\R^{N_p}\mbox{ such that:}}\\
\nu\lambda\bbA\,\ul{U}^n-\lambda\bbB^T\ul{P}^n&=&\lambda\ul{F}_u\\
\lambda\bbB\,\ul{U}^n+\nu^{-1}\bbM\,\ul{P}^n&=&\nu^{-1}\bbM\,\ul{P}^{n-1}
\end{array}
\right..
\end{equation}
In this situation, one needs to prove that the iterative solver is actually converging, and to provide a stopping criterion.
\begin{rmk} 
The matrix of algorithm \eqref{eq:StokesFV2h-Iterative} is similar to the one of the first-order artificial compressibility algorithm \cite{Shen95,GuMi15,GuMi17}. In \cite{GuMi15} and \cite{GuMi17}, the space discretization is done with the MAC approximation on a Cartesian grid \cite{HaWe65}. In \cite{GuMi17}, it is also done with the Taylor-Hood finite element \cite{TaHo73}.
\end{rmk}
\begin{thm}[Convergence] \label{thm-cv_algo}
Assume that the matrices $\bbA\in\R^{N_u\times N_u}$ and $\bbM\in\R^{N_p\times N_p}$ are symmetric positive-definite. Then, the sequence $(\ul{U}^n,\ul{P}^n)_n$ of solutions to \eqref{eq:StokesFV2h-Iterative} is converging, and it converges to $(\ul{U}^\infty,\ul{P}^\infty)\in\R^{N_u}\times\R^{N_p}$ which is governed by
\begin{equation}\label{eq:StokesP1P0_FVh}
\left\{
\begin{array}{rcl}
\multicolumn{3}{l}{\mbox{Find }(\ul{U}^\infty,\ul{P}^\infty)\in\R^{N_u}\times\R^{N_p}\mbox{ such that:}}\\
\nu\bbA\,\ul{U}^\infty-\bbB^T\ul{P}^\infty&=&\ul{F}_u\\
\bbB\,\ul{U}^\infty&=&0\\
\multicolumn{3}{l}{\bbM\,(\ul{P}^\infty-\ul{P}^{-1})\in \im(\bbB)}
\end{array}.
\right.
\end{equation}
\end{thm}
\begin{proof} 
In the proof, we use the change of variable $\ul{Y}:=\bbM^\frac12\ul{P}\in\R^{N_p}$, and specifically $\ul{Y}^n = \bbM^\frac12 \ul{P}^n$ for $n\ge-1$. Then let $\delta\ul{P}^n=\ul{P}^n-\ul{P}^{n-1}$ and $\delta\ul{Y}^n:=(\bbM^\frac12\delta \ul{P}^n)$ for $n\geq 0$.
 \\
1/ We notice first that the following recursive relation holds:
\begin{equation}\label{eq:deltaP}
\forall n\geq 1,\quad(\bbM+\lambda\bbB\bbA^{-1}\bbB^T)\delta\ul{P}^n=\bbM\delta\ul{P}^{n-1},
\end{equation}
where the matrix $\lambda\bbB\bbA^{-1}\bbB^T$ is a symmetric positive matrix. By construction, $\delta\ul{Y}^n = -\nu\lambda\bbM^{-\frac12}\bbB\,\ul{U}^n \in \im(\bbM^{-\frac12}\bbB)$. Since 
$\R^{N_p}=\im(\bbM^{-\frac12}\bbB)\stackrel{\perp}{\oplus}\ker(\bbB^T\bbM^{-\frac12})$, one has $\delta\ul{Y}^n\in\ker(\bbB^T\bbM^{-\frac12})^\perp$ and \eqref{eq:deltaP} reads:
\begin{equation}\label{eq:Yn}
\forall n\geq1,\quad(\mathbb{I}+\lambda\bbM^{-\frac12}\bbB\bbA^{-1}\bbB^T\bbM^{-\frac12})\delta\ul{Y}^n=\delta\ul{Y}^{n-1}.
\end{equation}
The iteration matrix $(\mathbb{I}+\lambda\bbM^{-\frac12}\bbB\bbA^{-1}\bbB^T\bbM^{-\frac12})^{-1}$ restricted to $\ker(\bbB^T\bbM^{-\frac12})^\perp$ is a symmetric positive definite matrix, whose eigenvalues are all strictly smaller than $1$. Hence, its spectral radius $\rho_\lambda$ is strictly smaller than $1$. \\
2/ We deduce that it exists $C_\lambda\in(0,1)$ independent of the data and of $n$ such that for all $n\geq 1$, $|\delta\ul{Y}^n|_2\leq C_\lambda|\delta\ul{Y}^{n-1}|_2$. Thus, the series $\sum_{n\geq 1}\delta\ul{Y}^n$ is a convergent series. Noticing that $\sum_{n=0}^N\delta\ul{Y}^n=\ul{Y}^N-\ul{Y}^{-1}$, we conclude that the sequence $(\ul{Y}^n)_n$ is a convergent sequence and we set $\ul{Y}^\infty=\lim_{n\rightarrow\infty}\ul{Y}^n$. By construction $\ul{Y}^\infty-\ul{Y}^{-1}\in\im(\bbM^{-\frac12}\bbB)=\ker(\bbB^T\bbM^{-\frac12})^\perp$, and for all $n$, the algorithm doesn't change the orthogonal projection of $\ul{Y}^n$ onto $\ker(\bbB^T\bbM^{-\frac12})$, which is equal to that of $\ul{Y}^{-1}$ onto $\ker(\bbB^T\bbM^{-\frac12})$.\\
3/ Because $\bbM^\frac12$ is invertible, we obtain that the sequence $(\ul{P}^n)_n$ is also a convergent sequence and we call $\ul{P}^\infty=\bbM^{-\frac12}\ul{Y}^\infty$ its limit. In the same way, noting that $\bbA$ is invertible, we obtain that the sequence $(\ul{U}^n)_n$ is a convergent sequence and we call $\ul{U}^\infty=\nu^{-1}\bbA^{-1}(\ul{F}_u+\bbB^T\ul{P}^\infty)$ its limit. \\
Passing to the limit, it holds that $\bbB\,\ul{U}^\infty=0$ and $\bbM(\ul{P}^\infty-\ul{P}^{-1})\in \im(\bbB)$.
\end{proof}
\begin{rmk} In the above proof, we notice that the mapping $\lambda\mapsto\rho_\lambda$ is monotonically decreasing. In section \ref{sec:resu} dedicated to the numerical experiments, we set the value of $\lambda$ equal to $1$ or $10$.
\end{rmk}
\begin{rmk} There are two exclusive cases for Problem \eqref{eq:StokesP1P0_FVh}:
\begin{itemize}
\item either, $\im(\bbB)=\R^{N_p}$: a discrete inf-sup condition holds. Uniqueness of $\ul{P}^n$ is guaranteed and the last line is trivial\,; 
\item or, $\im(\bbB)\subsetneq\R^{N_p}$: there is no discrete inf-sup condition. However in this case, the last line guarantees the uniqueness of $\ul{P}^n$. As a matter of fact, the algorithm does not change the orthogonal projection of $\ul{P}^n$ onto $\ker(\bbB^T)$, which is equal to that of $\ul{P}^{-1}$ onto $\ker(\bbB^T)$, where orthogonality is understood with respect to the inner product $(\bbM\cdot|\cdot)$.
\end{itemize}
\end{rmk}
Let us consider that the assumptions of theorem~\ref{thm-cv_algo} hold. 
A critical ingredient is to define the stopping criterion, in particular what are the relevant quantities of interest. For that, we use the expression $\bbM\,\delta\ul{P}^n=-\nu\lambda\bbB\ul{U}^n$. We recall that $\bbB\ul{U}^n$ stands for the divergence of the discrete velocity, and that one aims at diminishing the value of its norm. Considering again the auxiliary variable $\delta\ul{Y}^n$ introduced in the proof of theorem~\ref{thm-cv_algo}, we find:
\[ \lambda^2|\bbM^{-1/2}\bbB\ul{U}^n|_2^2 = \nu^{-2}|\delta\ul{Y}^n|_2^2. \]
As a consequence of the proof of the previous theorem, we infer that the sequence of norms $(|\bbM^{-1/2}\bbB\ul{U}^n|_2)_n$ is strictly monotonically decreasing.  \\
Hence, denoting by $|\cdot|_{\bbA}:\ul{U}\mapsto (\bbA\ul{U}|\ul{U})^{1/2}$ the norm of the discrete velocity, resp. $|\cdot|_{\bbM}:\ul{P}\mapsto (\bbM\ul{P}|\ul{P})^{1/2}$ the norm of the discrete pressure, we define the stopping criterion by comparing the adimensionalized quantities $|\delta\ul{P}^n|_{\bbM}$ and $|\ul{U}^n|_{\bbA}$ (we note that both quantities are easily computable). This amounts to setting the stopping criterion as
\begin{equation}\label{stopping_criterion}
|\delta\ul{P}^n|_{\bbM} \le \epsilon \, |\ul{U}^n|_{\bbA}
 \end{equation}
for ad hoc $\epsilon >0$.  In our numerical experiments, we fixed $n=1$ or $n=8$ with $\epsilon=10^{-12}$.
\section{Discretization and stability estimates}\label{sec:ananum}
\textcolor{cyan}{We solve numerically the new variational formulation, written like in Problem \eqref{eq:StokesFV2}. Below, we consider homogeneous Dirichlet boundary conditions (see \ref{sec:CLDNH} for nonhomogeneous Dirichlet boundary conditions)}.
\subsection{Discretizations}
Consider $(\Tcal_h)_h$ a simplicial triangulation sequence of $\Om$. For all $D\subset\R^d$, and $k\in\N$, we call $P^k(D)$ the set of order $k$ polynomials on $D$, $\bP^k(D)=(P^k(D))^d$, and we consider the broken polynomial space:
\[
P^k_{disc}(\Tcal_h):=\left\{q\in L^2(\Om);\quad \forall T\in\Tcal_h,\,q_{|T}\in P^k(T)\right\}.
\]
We use the notation $P^0(\Tcal_h)$ for $P^0_{disc}(\Tcal_h)$. 
We define the \textcolor{cyan}{conforming} spaces of $P^k$-Lagrange functions:
\[
\begin{array}{rcl}
V_h^k&:=&\ds\left\{v_h\in H^1(\Om);\quad\forall T\in\Tcal_h,\,v_{h|T}\in P^k(T)\right\}\\
 V_{0,h}^k&:=&\ds\left\{v_h\in V_h^k;\quad v_{h|\pa\Om}=0\right\}.
\end{array}
\]
We set $\bV_h^k:=(V_h^k)^d$ and $\bV_{0,h}^k:=(V_{0,h}^k)^d$. We call $Q_h^k:= P_{disc}^k(\Tcal_h) \cap L^2_{zmv}(\Om)$. 
Notice that $\bV_{0,h}^k\times Q_h^{k'}\subset\bH^1_0(\Om)\times L^2_{zmv}(\Om)$ and $\div\bV_{0,h}^k\subset Q_h^{k-1}$.

Let $k\ge1$. The (conforming) discretization of Problem \eqref{eq:StokesFV2} with \textcolor{cyan}{the Taylor-Hood} $\bP^k\times P^{k-1}_{disc}$ finite element \cite{TaHo73} reads:
\begin{equation}
\label{eq:StokesFV2h}
\left\{\begin{array}{l}
\mbox{Find }(\uvec_h,p_h)\in\bV_{0,h}^k\times Q_h^{k-1}\mbox{ s.t. for all }(\vvec_h,q_h)\in\bV_{0,h}^k\times Q_h^{k-1}:
\\
\nu\,\lambda(\uvec_h,\vvec_h)_{\bH^1_0(\Om)}-\lambda(p_h,\div\vvec_h)_{L^2(\Om)} = \lambda\langle\fvec,\vvec_h\rangle_{\textcolor{cyan}{\bH^{-1}(\Om),\bH^1_0(\Om)}},\\
\lambda(q_h,\div\uvec_h)_{L^2(\Om)} + \nu^{-1}\,(p_h,q_h)_{L^2(\Om)} = \nu^{-1}\,(z_{\fvec},q_h)_{L^2(\Om)}.
\end{array}\right.
\end{equation}
Writing $\fvec=-\nu\Delta\uvec+\grad z_{\fvec}$ (recall (\ref{eq:Stokes}) with $z_\fvec=p$), Problem \eqref{eq:StokesFV2h} reads:
\begin{equation}
\label{eq:StokesFV2h-bis}
\left\{\begin{array}{l}
\mbox{Find }(\uvec_h,p_h)\in\bV_{0,h}^k\times Q_h^{k-1}\mbox{ s.t. for all }(\vvec_h,q_h)\in\bV_{0,h}^k\times Q_h^{k-1}:
\\
\nu\,\lambda(\uvec_h,\vvec_h)_{\bH^1_0(\Om)}-\lambda(p_h,\div\vvec_h)_{L^2(\Om)} = \nu\,\lambda(\uvec,\vvec_h)_{\bH^1_0(\Om)}\\
\hskip 65truemm-\lambda(z_{\fvec},\div\vvec_h)_{L^2(\Om)},\\
\lambda(q_h,\div\uvec_h)_{L^2(\Om)} + \nu^{-1}\,(p_h,q_h)_{L^2(\Om)} =\nu^{-1} (z_{\fvec},q_h)_{L^2(\Om)}.
\end{array}\right.
\end{equation}
\begin{rmk}
It is well known that discretizing \eqref{eq:Stokes-FV} with the \textcolor{cyan}{Taylor-Hood $\bP^k\times P^{k-1}_{disc}$ finite element}, $k\geq 1$, is not stable on all shape regular meshes \cite{BoBF13}. A wide range of strategies to get round this problem have been explored for years. Below is a quick overview of these strategies:
\begin{itemize}
\item The discrete velocity space can be enriched see \cite{BeRa81,ArBF84,BaVa02,CGRT18,LiZi22} for $k=1$.
\item The discrete variational formulation can be stabilized, see \cite{BoDG06,BaVa10,ZhTa15,AlGR18} for $k=1$, \cite{AABR13} for $k=2$.%
\item The mesh can be designed in such a way that the {uniform} discrete inf-sup condition applies, see \cite{FGNZ22} for $k=1$, \cite{ScVo85,Zhan05,BBCG23} for $k>1$.
\end{itemize}
\end{rmk}

\subsection{Error estimate when $z_{\fvec}$ is known}\label{ss:cv-z-known}
Let $\Pi_{h,cg}$ be the $\bL^2(\Om)$-orthogonal projection operator from $\bL^2(\Om)$ to $\bV_{0,h}^k$ and $\Pi_{h,dg}$ be the $L^2(\Om)$-orthogonal projection operator from $L^2(\Om)$ to $P^{k-1}_{disc}(\Tcal_h)$. We have the
\begin{prop}
\label{pro:pressure-robust-estimate}
Let $(\uvec,p)\in\Xcal$ be the exact solution and $(\uvec_h,p_h)\in\Xcal_h$ be the solution to Problem \eqref{eq:StokesFV2h}. We have the following error estimate:
\begin{equation}
\label{eq:pressure-robust-estimate}
\|(\uvec_h-\uvec,p_h-\Pi_{h,dg} z_{\fvec})\|_{\Xcal,{\sqrt{\lambda}}\nu}
\leq \sqrt{1+d\lambda}\,\|\uvec-\Pi_{h,cg}\uvec\|_{\bH^1_0(\Om)}.
\end{equation}
\textcolor{cyan}{Suppose that $\uvec\in (P^k(\Omega))^d$. Then $\uvec_h=\uvec$ and $p_h= \Pi_{h,dg}z_{\fvec}$}.
\end{prop}
Notice that the error estimate for the velocity is {independent} of the pressure. Hence, the (conforming) discretization of Problem \eqref{eq:StokesFV2} with $\bP^k\times P^{k-1}_{disc}$ \textcolor{cyan}{and for which $z_{\fvec}$ is known} is a so called {\it pressure robust method} for which the discrete velocity $\uvec_h$ tends to $\uvec$ independently of $\nu$. The discrete pressure $p_h$ tends to $\Pi_{h,dg}z_{\fvec}$ all the faster the smaller $\nu$ is. Let us now prove Proposition \ref{pro:pressure-robust-estimate}.
\begin{proof}
Setting $z_{\fvec,h}=\Pi_{h,dg}z_{\fvec}$, we have: $(z_{\fvec},\div\vvec_h)_{L^2(\Om)}=(z_{\fvec,h},\div\vvec_h)_{L^2(\Om)}$ and $(z_{\fvec},q_h)_{L^2(\Om)}=(z_{\fvec,h},q_h)_{L^2(\Om)}$ for all $\vvec_h\in\bV_{0,h}^k$, and all $q_h\in Q_h^{k-1}$.
Summing both equations of Problem \eqref{eq:StokesFV2h-bis} and reshuffling terms, it comes:
\begin{equation}
\label{eq:StokesFV2h-bis-sum}
\left\{\begin{array}{l}
\mbox{Find }(\uvec_h,p_h)\in\bV_{0,h}^k\times Q_h^{k-1}\mbox{ s.t. for all }(\vvec_h,q_h)\in\bV_{0,h}^k\times Q_h^{k-1}: 
\\
\nu\,\lambda(\uvec_h-\uvec,\vvec_h)_{\bH^1_0(\Om)}+\nu^{-1}\,(p_h-z_{\fvec,h},q_h)_{L^2(\Om)}\\ \quad=\lambda(p_h-z_{\fvec,h},\div\vvec_h)_{L^2(\Om)}-\lambda(q_h,\div\uvec_h)_{L^2(\Om)}.
\end{array}\right.
\end{equation}
Choosing $\vvec_h=\uvec_h-\Pi_{h,cg}\uvec=(\uvec_h-\uvec)+(\uvec-\Pi_{h,cg}\uvec)$ and $q_h=p_h-z_{\fvec,h}$ in \eqref{eq:StokesFV2h-bis-sum}, and dividing by $\lambda\nu$, we obtain:
\[
\begin{array}{c}
\|\uvec_h-\uvec\|_{\bH^1_0(\Om)}^2+(\,\uvec_h-\uvec,\uvec
-\Pi_{h,cg}\uvec\,)_{\bH^1_0(\Om)}
+\lambda^{-1}\nu^{-2}\|p_h-z_{\fvec,h}\|_{L^2(\Om)}^2\\
=-\nu^{-1}(p_h-z_{\fvec,h},\div\Pi_{h,cg}\uvec)_{L^2(\Om)}.
\end{array}
\]
Hence, using Cauchy-Schwarz inequality, we get:
\[
\begin{array}{l}
\|\uvec_h-\uvec\|_{\bH^1_0(\Om)}^2+\lambda^{-1}\nu^{-2}\|p_h-z_{\fvec,h}\|_{L^2(\Om)}^2\\
\leq\|\uvec_h-\uvec\|_{\bH^1_0(\Om)}\,\|\uvec-\Pi_{h,cg}\uvec\|_{\bH^1_0(\Om)}
+\,\nu^{-1}\|p_h-z_{\fvec,h}\|_{L^2(\Om)}\,\|\div\Pi_{h,cg}\uvec\|_{L^2(\Om)} \\
\le {(\|\uvec_h-\uvec\|_{\bH^1_0(\Om)}^2+\lambda^{-1}\nu^{-2}\|p_h-z_{\fvec,h}\|_{L^2(\Om)}^2)^{1/2} \times} \\
 \hskip 30truemm {(\|\uvec-\Pi_{h,cg}\uvec\|_{\bH^1_0(\Om)}^2 + \lambda\|\div\Pi_{h,cg}\uvec\|_{L^2(\Om)}^2 )^{1/2}}
\end{array}
\]
Finally $\|\div\Pi_{h,cg}\uvec\|_{L^2(\Om)}\leq \sqrt{d}\,\|\uvec-\Pi_{h,cg}\uvec\|_{\bH^1_0(\Om)}$, so we have:
\[
\|(\uvec_h-\uvec,p_h-z_{\fvec,h})\|_{\Xcal,{\sqrt{\lambda}}\nu}
\leq{\sqrt{1+d\lambda}}\,\|\uvec-\Pi_{h,cg}\uvec\|_{\bH^1_0(\Om)}.
\]
\textcolor{cyan}{If $\uvec\in (P^k(\Omega))^d$, then $\Pi_{h,cg}\uvec=\uvec$ and the right-hand side 
is equal to $0$. Hence, $\uvec_h=\uvec$ and $p_h= z_{\fvec,h}$.}
\end{proof}
\subsection{Error estimate when $z_{\fvec}$ is not known}
If $z_{\fvec}$ is not known explicitly, let us assume that we have at hand an approximation ${p_h^{old}}\in Q_h^{k-1}$ and use it to solve the following approximation of Problem \eqref{eq:StokesFV2}:
\begin{equation}
\label{eq:StokesFV3}
\left\{\begin{array}{l}
\mbox{Find $(\tilde{\uvec},\tilde{p})\in\Xcal$ s.t. for all }(\vvec,q)\in\Xcal:\\
\nu\,\lambda(\tilde{\uvec},\vvec)_{\bH^1_0(\Om)}-\lambda(\tilde{p},\div\vvec)_{L^2(\Om)} =\lambda\langle\fvec,\vvec\rangle_{\textcolor{cyan}{\bH^{-1}(\Om),\bH^1_0(\Om)}} ,\\
\lambda(q,\div\tilde{\uvec})_{L^2(\Om)} + \nu^{-1}\,(\tilde{p},q)_{L^2(\Om)} = \nu^{-1}\,({p_h^{old}},q)_{L^2(\Om)}.
\end{array}\right.
\end{equation}
Let us write again $\fvec=-\nu\Delta\uvec+\grad z_{\fvec}$, so that Problem \eqref{eq:StokesFV3} reads:
\begin{equation}
\label{eq:StokesFV3-bis}
\left\{\begin{array}{l}
\mbox{Find $(\tilde{\uvec},\tilde{p})\in\Xcal$ s.t. for all }(\vvec,q)\in\Xcal:\\
\nu\,\lambda(\tilde{\uvec}-\uvec,\vvec)_{\bH^1_0(\Om)}-\lambda(\tilde{p}-z_{\fvec},\div\vvec)_{L^2(\Om)} =  0,\\
\lambda(q,\div\tilde{\uvec})_{L^2(\Om)} + \nu^{-1}\,(\tilde{p}-{p_h^{old}},q)_{L^2(\Om)}=0.
\end{array}\right.
\end{equation}
\begin{prop}\label{pro:zfh}
Let $(\uvec,p)\in\Xcal$ be the exact solution and $(\tilde{\uvec},\tilde{p})\in\Xcal$ be the solution to Problem \eqref{eq:StokesFV3}. We have the following estimates:
\begin{equation}
\label{eq:estim-p-ptilde}
\left\{
\begin{array}{rcl}
\|\tilde{p}-z_{\fvec}\|_{L^2(\Om)}&\leq&\|z_{\fvec}-{p_h^{old}}\|_{L^2(\Om)},\\
\|\tilde{\uvec}-\uvec\|_{\bH^1_0(\Om)}&\leq&{(\sqrt{\lambda}\nu)^{-1}}\|z_{\fvec}-{p_h^{old}}\|_{L^2(\Om)}.
\end{array}
\right.
\end{equation}
\end{prop}
\begin{proof}
Choosing $\vvec=\tilde{\uvec}-\uvec$ and $q=\tilde{p}-z_{\fvec}$ in \eqref{eq:StokesFV3-bis}, summing both equations and dividing by $\lambda\nu$, it comes:
\[
\|\tilde{\uvec}-\uvec\|_{\bH^1_0(\Om)}^2+\lambda^{-1}\nu^{-2}\,(\tilde{p}-{p_h^{old}},\tilde{p}-z_{\fvec})_{L^2(\Om)} =0.
\]
Writing $\tilde{p}-{p_h^{old}}=(\tilde{p}-z_{\fvec})+(z_{\fvec}-{p_h^{old}})$, we obtain:
\[
\begin{array}{r}
\|(\tilde{\uvec}-\uvec,\tilde{p}-z_{\fvec})\|_{\Xcal,{\sqrt{\lambda}}\nu}^2=-\lambda^{-1}\nu^{-2}(z_{\fvec}-{p_h^{old}},\tilde{p}-z_{\fvec})_{L^2(\Om)}
\\
\leq\lambda^{-1}\nu^{-2}\|z_{\fvec}-{p_h^{old}}\|_{L^2(\Om)}\,\|\tilde{p}-z_{\fvec}\|_{L^2(\Om)}.
\end{array}
\]
We obtain successively both estimates \eqref{eq:estim-p-ptilde}.
\end{proof}
The discretization of Problem \eqref{eq:StokesFV3} with $\bP^k\times P_{disc}^{k-1}$ finite element reads:
\begin{equation}
\label{eq:StokesFV3h}
\left\{\begin{array}{l}
\mbox{Find $(\tilde{\uvec}_h,{p_h^{new}})\in\bV_{0,h}^k\times Q_h^{k-1}$ s.t. for all }(\vvec_h,q_h)\in\bV_{0,h}^k\times Q_h^{k-1}:\\
 \nu\,\lambda(\tilde{\uvec}_h,\vvec_h)_{\bH^1_0(\Om)}-\lambda({p_h^{new}},\div\vvec_h)_{L^2(\Om)} =\lambda\langle \fvec,\vvec_h\rangle_{\textcolor{cyan}{\bH^{-1}(\Om),\bH^1_0(\Om)}},\\
\lambda(q_h,\div\tilde{\uvec}_h)_{L^2(\Om)} + \nu^{-1}\,({p_h^{new}},q_h)_{L^2(\Om)} = \nu^{-1}\,({p_h^{old}},q_h)_{L^2(\Om)}.
\end{array}\right.
\end{equation}
\begin{thm}\label{th:zfh}
Let $(\uvec,p)\in\Xcal$ be the exact solution and $(\tilde{\uvec}_h,{p_h^{new}})\in\Xcal$ be the solution to Problem \eqref{eq:StokesFV3h}. We have the following error  estimate:
\begin{equation}
\label{eq:estim-p-ptildeh}
\begin{array}{r}
\|(\tilde{\uvec}_h-\uvec,{p_h^{new}}- \Pi_{h,dg}z_{\fvec})\|_{\Xcal,{\sqrt{\lambda}}\nu}\leq \sqrt{1+d\lambda}\,\|\uvec-\Pi_{h,cg}\uvec\|_{\bH^1_0(\Om)}\\
+{(\sqrt{\lambda}\nu)^{-1}}\,\| \Pi_{h,dg}z_{\fvec}-{p_h^{old}}\|_{L^2(\Om)}.
\end{array}
\end{equation}
\textcolor{cyan}{Suppose that $\uvec\in (P^k(\Omega))^d$. Then $\Pi_{h,cg}\uvec=\uvec$ and we obtain the estimate below:}
\begin{equation}\label{eq:u=0-errorestimate}
\|(\tilde{\uvec}_h-\uvec,{p_h^{new}}- \Pi_{h,dg}z_{\fvec})\|_{\Xcal,{\sqrt{\lambda}}\nu}\leq
{(\sqrt{\lambda}\nu)^{-1}}\,\|\Pi_{h,dg}z_{\fvec}-{p_h^{old}}\|_{L^2(\Om)}.
\end{equation}
\end{thm}
In particular, if ${p_h^{old}}$ is a {good} approximation of $\Pi_{h,dg}z_{\fvec}$, the solution $(\tilde{\uvec}_h,{p_h^{new}})$ {is also a good approximation of} $(\uvec,\Pi_{h,dg}z_{\fvec})$. Interestingly, the above with $k=1$ corresponds to the $\bP^1\times P^0$ finite element pair, which is known to be unstable for the discretization of the usual variational formulation (\ref{eq:Stokes-FV}) of the Stokes problem. Let us prove Theorem \ref{th:zfh}.
\begin{proof}
Setting $z_{\fvec,h}=\Pi_{h,dg}z_{\fvec}$, we have: $\langle \fvec,\vvec_h\rangle_{\textcolor{cyan}{\bH^{-1}(\Om),\bH^1_0(\Om)}}=\nu\,(\uvec,\vvec_h)_{\bH^1_0(\Om)}-(z_{\fvec,h},\div\vvec_h)_{L^2(\Om)}$ for all $\vvec_h\in\bV_{0,h}^k$. Hence, Problem \eqref{eq:StokesFV3h} can be written as:
 \begin{equation}
\label{eq:StokesFV3h-bis}
\left\{\begin{array}{l}
\mbox{Find $(\tilde{\uvec}_h,{p_h^{new}})\in\bV_{0,h}^k\times Q_h^{k-1}$ s.t. for all }(\vvec_h,q_h)\in\bV_{0,h}^k\times Q_h^{k-1}:\\
\nu\,\lambda(\tilde{\uvec}_h-\uvec,\vvec_h)_{\bH^1_0(\Om)}-\lambda({p_h^{new}}-z_{\fvec,h},\div\vvec_h)_{L^2(\Om)}=0,\\
\lambda(q_h,\div\tilde{\uvec}_h)_{L^2(\Om)} + \nu^{-1}\,({p_h^{new}}-{p_h^{old}},q_h)_{L^2(\Om)}=0.
\end{array}\right.
\end{equation}
Choosing $(\vvec_h,q_h)=(\tilde{\uvec}_h-\Pi_{h,cg}\uvec,{p_h^{new}}-z_{\fvec,h})$  in  Problem \eqref{eq:StokesFV3h-bis}, we have:
\[
\left\{\begin{array}{l}
\nu\,\lambda(\tilde{\uvec}_h-\uvec,\tilde{\uvec}_h-\Pi_{h,cg}\uvec)_{\bH^1_0(\Om)}-\lambda({p_h^{new}}-z_{\fvec,h},\div(\tilde{\uvec}_h-\Pi_{h,cg}\uvec))_{L^2(\Om)} = 0,\\
\lambda({p_h^{new}}-z_{\fvec,h},\div\tilde{\uvec}_h)_{L^2(\Om)} + \nu^{-1}\,({p_h^{new}}-{p_h^{old}},{p_h^{new}}-z_{\fvec,h})_{L^2(\Om)} = 0.
\end{array}\right.
\]
Summing both equations and dividing by $\lambda\nu$, it now comes:
\[
\begin{array}{c}
(\tilde{\uvec}_h-\uvec,\tilde{\uvec}_h-\Pi_{h,cg}\uvec)_{\bH^1_0(\Om)} + \lambda^{-1}\nu^{-2}\,({p_h^{new}}-{p_h^{old}},{p_h^{new}}-z_{\fvec,h})_{L^2(\Om)} \\
+\nu^{-1}({p_h^{new}}-z_{\fvec,h},\div\Pi_{h,cg}\uvec)_{L^2(\Om)}=0.
\end{array}
\]
Noticing that $\tilde{\uvec}_h-\Pi_{h,cg}\uvec=(\tilde{\uvec}_h-\uvec)+(\uvec-\Pi_{h,cg}\uvec)$ and ${p_h^{new}}-{p_h^{old}}=({p_h^{new}}-z_{\fvec,h})+(z_{\fvec,h}-{p_h^{old}})$, we get:
\[
\begin{array}{c}
\|(\tilde{\uvec}_h-\uvec,{p_h^{new}}-z_{\fvec,h})\|_{\Xcal,{\sqrt{\lambda}}\nu}^2+(\tilde{\uvec}_h-\uvec,\uvec-\Pi_{h,cg}\uvec)_{\bH^1_0(\Om)} \\
+\lambda^{-1}\nu^{-2}\,(z_{\fvec,h}-{p_h^{old}},{p_h^{new}}-z_{\fvec,h})_{L^2(\Om)}
+\nu^{-1}\,({p_h^{new}}-z_{\fvec,h},\div\Pi_{h,cg}\uvec)_{L^2(\Om)}=0.
\end{array}
\]
Using Cauchy-Schwarz inequality, we deduce that:
\[
\begin{array}{c}
\|(\tilde{\uvec}_h-\uvec,{p_h^{new}}-z_{\fvec,h})\|_{\Xcal,{\sqrt{\lambda}}\nu}^2\leq\|\tilde{\uvec}_h-\uvec\|_{\bH^1_0(\Om)}\,\|\uvec-\Pi_{h,cg}\uvec\|_{\bH^1_0(\Om)} \\
+{(\sqrt{\lambda}\nu)^{-1}}\,\|{p_h^{new}}-z_{\fvec,h}\|_{L^2(\Om)}\left({(\sqrt{\lambda}\nu)^{-1}}\,\|z_{\fvec,h}-{p_h^{old}}\|_{L^2(\Om)}+{\sqrt{\lambda}}\|\div\Pi_{h,cg}\uvec\|_{L^2(\Om)}\right).
\end{array}
\]
Using again Cauchy-Schwarz inequality as in the proof of Proposition~\ref{pro:pressure-robust-estimate}, together with $\|\div\Pi_{h,cg}\uvec\|_{L^2(\Om)}\leq \sqrt{d}\,\|\uvec-\Pi_{h,cg}\uvec\|_{\bH^1_0(\Om)}$, we obtain the estimates \eqref{eq:estim-p-ptildeh} and \eqref{eq:u=0-errorestimate}.
\end{proof}
\subsection{\textcolor{cyan}{Discussion on the error indicators}}
{As we noted after stating Theorem~\ref{th:zfh}, if the error indicator on the pressure $\|z_{\fvec,h}-{p_h^{old}}\|_{L^2(\Om)}$ is small, so is $\|z_{\fvec,h}-{p_h^{new}}\|_{L^2(\Om)}$. Indeed, assuming that the errors on the velocities are negligible compared to these indicators, one derives from (\ref{eq:estim-p-ptildeh}) the bound $\|z_{\fvec,h}-{p_h^{new}}\|_{L^2(\Om)} \le \|z_{\fvec,h}-{p_h^{old}}\|_{L^2(\Om)}$. Interestingly, this situation is likely to occur when the viscosity is small, since there is a multiplicative factor equal to ${(\sqrt{\lambda}\nu)^{-1}}$ in front of the error indicators on the pressure. This further justifies the use of the iterative algorithm presented in Section~\ref{sec:algo} in this case.}

\section{Numerical results}\label{sec:resu}
\textcolor{cyan}{In this section, the computations are performed on a Dell Precision 3581 Intel Core i7 laptop.} We propose some numerical experiments \textcolor{cyan}{with $g=0$}, depending on whether or not $z_{\fvec}$ is known explicitly. In the latter case, we first compute some approximation ${p_h^{old}}$. In principle, either a classical, conforming or nonconforming, discretization can be used to compute ${p_h^{old}}$. \textcolor{cyan}{Hence, we apply an iterative approach with a nonconforming initial computation followed by (several iterations of) our new formulation. we call the second part the post-processing. For each test-case, $(\uvec,p)$ is given, whereas the viscosity $\nu$ can vary (see the discussion after (\ref{eq:normX})).} The numerical results are obtained on a github platform, 
implemented in Octave language, see \cite{StokesExplicitTC}.
\subsection{Resolution algorithm} \textcolor{cyan}{Consider nonconforming discretization of the classical variational formulation (\ref{eq:Stokes-FV}), cf. \cite{CrRa73}. Let $\bV_{0,h}^{nc}\not\subset\bH^1_0(\Om)$) be the discrete velocity space, and $Q_h$ be the discrete pressure space. Let $(\psi_i)_{i=1}^{N_u}$ be a basis of $\bV_{0,h}^{nc}$, and $(\phi_i)_{i=1}^{N_p}$ be a basis of $Q_h$. We set: $\ul{U}:=(\ul{U}_i)_{i=1}^{N_u}$ where $\uvec_h:=\sum_{i=1}^{N_u}\ul{U}_i\psi_i$, $\ul{P}:=(\ul{P}_i)_{i=1}^{N_p}$ where $p_h:=\sum_{i=1}^{N_p}\ul{P}_i\phi_i $, and $\ul{F}_u:=(\ul{F}_{u,i})_{i=1}^{N_u}$. We let $\ell_\fvec\in\Lcal(\bV_{0,h}^{nc},\R)$ be such that $\forall\vvec_h\in\bV_{0,h}^{nc}$, $\ell_\fvec(\vvec_h)=(\fvec,\vvec_h)_{\bL^2(\Om)}$ if $\fvec\in\bL^2(\Om)$, $\ell_\fvec(\vvec_h)=\langle\fvec,\Jcal_h(\vvec_h)\rangle_{\textcolor{cyan}{\bH^{-1}(\Om),\bH^1_0(\Om)}}$ if $\fvec\not\in\bL^2(\Om)$, where $\Jcal_h:\bV_{0,h}^{nc}\rightarrow\bV_{0,h}^k$ is for instance an averaging operator \cite[\S 22.4.1]{ErGu21-I}. Then $\ul{F}_{u,i}=\ell_\fvec(\psi_i)$.}
Let $\bbA\in\R^{N_u}\times\R^{N_u}$ be the velocity stiffness matrix, $\bbB\in\R^{N_p}\times\R^{N_u}$ be the velocity-pressure coupling matrix and $\bbM\in\R^{N_p}\times \R^{N_p}$ be the pressure mass matrix. The linear system to be solved is
\begin{equation}\label{eq:StokesClassique_FVh}
\left\{
\begin{array}{rcl}
\multicolumn{3}{l}{\mbox{Find }(\ul{U},\ul{P})\in\R^{N_u}\times\R^{N_p}\mbox{ such that:}}\\
\nu\bbA\,\ul{U}-\bbB^T\ul{P}&=&\ul{F}_u\\
\bbB\,\ul{U}&=&0
\end{array}.
\right.
\end{equation}
Let us set $\bbK=\bbB\,\bbA^{-1}\bbB^T\in\R^{N_p}\times\R^{N_p}$. The matrix $\bbK\in\R^{N_p}\times\R^{N_p}$ is a symmetric matrix, furthermore it is positive definite as soon as the kernel of $\bbB^T$ is reduced to~$\{0\}$. When it is the case, \textcolor{cyan}{the coupled velocity-pressure problem (\ref{eq:StokesClassique_FVh}) can be solved 
using the three + one steps below (the fourth step being straightforward):
\begin{equation}\label{eq:Uzawa}
\begin{array}{ll}
\mbox{Prediction:}&\mbox{Solve in }\ul{U}_{\star}\mbox{ such that }\nu\,\bbA\,\ul{U}_{\star}= \ul{F}_u.\\
\mbox{Pressure solver:}&\mbox{Solve in }\ul{P}\mbox{ such that }\bbK\,\ul{P}=\ul{F}_p-\nu\,\lambda\bbB\,\ul{U}_{\star}.\\
\mbox{Correction:}&\mbox{Solve in }\delta\ul{U}\mbox{ such that }\nu\,\bbA\,\delta\ul{U}=\bbB^T\,\ul{P}.
\\
\mbox{Update:}&\ul{U}=\delta\ul{U}+\ul{U}_{\star}.
\end{array}
\end{equation}
One can check easily that the above computed solution $(\ul{U},\ul{P})$ solves (\ref{eq:StokesClassique_FVh}). The pressure solver with matrix $\bbK$ is based on the Uzawa algorithm, which is the conjugate gradient algorithm in the context of the Stokes problem. It can be preconditioned by the inverse of the mass matrix associated to the discrete pressure $\bbM$ (see e.g. \cite[Lemma 5.9]{OlTy14}). Thanks to the uniform discrete inf-sup condition, the number of iterations of the conjugate gradient algorithm is independent of the meshsize. 
The matrices $\bbM$ and $\bbB$ are kept with all the $P_{disc}^{k-1}$ degrees of freedom for the discrete pressure. To take into account the zero mean value constraint, at each iteration of the preconditioned conjugate gradient algorithm, the discrete pressure is first computed in $L^2(\Om)$, then orthogonally projected in $L^2_{zmv}(\Om)$. We use the Cholesky factorization (computed once and for all) to solve linear systems with matrix $\bbA$.
}

Consider next the $\bP^k\times P^{k-1}_{disc}$ conforming discretization of the variational formulation (\ref{eq:StokesFV2h}) or (\ref{eq:StokesFV3h}). Let $(\psi_i)_{i=1}^{N_u}$ be the Lagrange basis of $\bV^k_{0,h}$ and $(\phi_i)_{i=1}^{N_p}$ be the basis of $Q^{k-1}_h$. We set: $\ul{U}:=(\ul{U}_i)_{i=1}^{N_u}$ where $\uvec_h:=\sum_{i=1}^{N_u}\ul{U}_i\psi_i$ and $\ul{P}:=(\ul{P}_i)_{i=1}^{N_p}$ where $p_h:=\sum_{i=1}^{N_p}\ul{P}_i\phi_i $. We set $\ul{F}_u:=(\ul{F}_{u,i})_{i=1}^{N_u}$ where $\ul{F}_{u,i}=\langle\fvec,\psi_i\rangle_{\textcolor{cyan}{\bH^{-1}(\Om),\bH^1_0(\Om)}}$ and $\ul{F}_p:=(\ul{F}_{p,i})_{i=1}^{N_p}$ where $\ul{F}_{p,i}=(z_\fvec,\phi_i)_{L^2(\Om)}$, cf. (\ref{eq:StokesFV2h}), or $\ul{F}_{p,i}=({p_h^{old}},\phi_i)_{L^2(\Om)}$, where the discrete pressure ${p_h^{old}}$ is an approximation of $z_\fvec$, cf. (\ref{eq:StokesFV3h}). The linear system to be solved is
\begin{equation}\label{eq:StokesFV2h-LS}
\left\{
\begin{array}{rcl}
\multicolumn{3}{l}{\mbox{Find }(\ul{U},\ul{P})\in\R^{N_u}\times\R^{N_p}\mbox{ such that:}}\\
\nu\lambda\bbA\,\ul{U}-\lambda\bbB^T\ul{P}&=&\lambda\ul{F}_u\\
\lambda\bbB\,\ul{U}+\nu^{-1}\bbM\,\ul{P}&=&\nu^{-1}\ul{F}_p
\end{array}
\right..
\end{equation}
\textcolor{cyan}{
In that case, we set $\bbA_{\lambda}=\bbA+\lambda\bbB^T\bbM^{-1}\bbB\in\R^{N_u}\times\R^{N_u}$, which is automatically a symmetric positive definite matrix. To solve the coupled velocity-pressure problem (\ref{eq:StokesFV2h-LS}), one relies on the two steps below :
\begin{equation}\label{eq:Uzawa-Uh}
\begin{array}{ll}
\mbox{Velocity solver:}&\mbox{Solve in }\ul{U}\mbox{ such that }\nu\,\bbA_{\lambda}\,\ul{U}= \ul{F}_u+\bbB^T\bbM^{-1}\ul{F}_p.\\
\mbox{Pressure solver:}&\mbox{Solve in }\ul{P}\mbox{ such that }\bbM\,\ul{P}=\ul{F}_p-\nu\lambda\bbB\,\ul{U}.
\end{array}
\end{equation}
One can check easily that the above computed solution $(\ul{U},\ul{P})$ solves (\ref{eq:StokesFV2h-LS}). The matrix $\bbA_\lambda$ can be assembled easily since $\bbM$ is a block-diagonal matrix. Using $\bP^1\times P^0$ discretization, the computation of $\ul{P}$ is straightforward. Again, matrices $\bbM$ and $\bbB$ are kept with all the $P_{disc}^{k-1}$ degrees of freedom for the discrete pressure. We proceed as before to take into account the zero mean value constraint.
}
\begin{rmk}\label{rem:cond}
We recall that solving the linear system (\ref{eq:StokesClassique_FVh}) via the solver (\ref{eq:Uzawa}) is not so straightforward, even with $\nu=1$. Solving the linear system (\ref{eq:StokesFV2h-LS}) is much easier.
\end{rmk}
\subsection{Settings}
We consider Problem \eqref{eq:Stokes} with homogeneous or non homogeneous boundary conditions in $\Om=(0,1)^2$. Let $(\uvec,p)$ be the exact solution, and $(\uvec_h,p_h)$ be the numerical solution. \textcolor{cyan}{We compare numerical methods, showing how the coercive $\bP^1\times P^0$ formulation can be used in a post-processing step (i.e. solving \eqref{eq:StokesFV3h} with $p_h^{old}$ known), improving then the approximation of the discrete velocity}.
On Figures \ref{fig:XYFlow} to \ref{fig:LowRegCPU6h}, and in Tables \ref{tab:XYFlow-CPU} to \ref{tab:LowReg-CPU} we give the following names to our numerical methods:
\begin{itemize}
\item {Method with Crouzeix-Raviart (CR)}: computations are made with the nonconforming Crouzeix-Raviart $\bP^1_{nc}\times P^0$ formulation \cite[Example 4]{CrRa73}, which is not a pressure robust method. We call $p_{nc}$ the resulting discrete pressure.
\item {Method with exact pressure (EP)}: computations are made with the coercive $\bP^1\times P^0$ formulation \eqref{eq:StokesFV2h} and $\textcolor{cyan}{\lambda=1}$, knowing the {exact pressure} $z_\fvec$ (i.e. we solve Problem \eqref{eq:StokesFV2h}). 
\item {Method with post-processing (Post)}: computations are made in {two steps}. In a first step, we approximate the pressure by the {CR-method}. Then, in a post processing step, we use this numerical pressure as the source term in the {EP-method} (i.e. we solve Problem \eqref{eq:StokesFV3h} with ${p_h^{old}}=p_{nc}$) and \textcolor{cyan}{$\lambda=1$ when $\nu=1$, $\lambda=10$ else. Unless the stopping criterion (\ref{stopping_criterion}) is achieved}, we iterate eight times the second step, updating ${p_h^{old}}$ at each new iteration. \textcolor{cyan}{The algorithm corresponds to Algorithm \eqref{eq:StokesFV2h-Iterative}.}
\end{itemize}
Let $N_T$ be the number of triangles. For the velocity, the number of unknowns is of order $N_u\approx 3\,N_T$  for the {CR-method} and $N_u\approx N_T$ for the {EP-method}. For the pressure, the number of unknowns is $N_p=N_T-1$ for both methods. As a consequence, there are roughly twice as many unknowns for the {CR-method} than there are for the {EP-method}. We report in Table~\ref{tab:NoU} the number of unknowns "$\#$ dof" for the numerical tests.
\begin{table}[h]
\begin{center}
\begin{tabular}{c||c|c||c|c}
$h$      & $\#$ dof CR & $\#$ dof EP & $\#$ dof CR & $\#$ dof EP \\
\hline
$1.00\times 10^{-1}$&  $1\,048$ &     $566$ &   $2\,184$&   $1\,114$ \\
$5.00\times 10^{-2}$&  $4\,376$ &  $2\,270$ &   $8\,064$&   $4\,074$ \\
$2.50\times 10^{-2}$& $17\,368$ &  $8\,846$ &  $30\,464$&  $15\,314$ \\
$1.25\times 10^{-2}$& $67\,816$ & $34\,230$ & $117\,664$&  $58\,994$\\
$6.25\times 10^{-3}$&$272\,624$ &$136\,954$ & $464\,448$& $232\,866$\\
\end{tabular}
\caption{Number of unknowns: first two test cases (left), last test case (right).}
\label{tab:NoU}
\end{center}
\end{table}

We propose three numerical examples based on manufactured solutions.
\textcolor{cyan}{Let $\Ical_{h,cg}$ be the interpolation operator from $\Ccal^0(\overline{\Omega})$ to $\bV_{0,h}^1$ and $\Ical_{h,nc}$ be the interpolation operator from $\Ccal^0(\overline{\Omega})$ to $\bV_{0,h}^{nc}$. We use the discrete $L^2$-error for the velocity and for the pressure:
\begin{equation}\label{eq:errUPL2}
\begin{array}{rcl}
\eps_0^\nu(\uvec_h)&:=&\|\Ical_h\uvec-\uvec_h\|_{\bL^2(\Om)}/\|(\uvec,p)\|_{\Xcal,\nu}\\
\eps_0^\nu(p_h)&:=&\nu^{-1}\|\Pi_{h,dg}p-p_h\|_{L^2(\Om)}/\|(\uvec,p)\|_{\Xcal,\nu}
\end{array},
\end{equation}
where $\Ical_h$ is either equal $\Ical_{h,nc}$ or $\Ical_{h,cg}$. 
For the EP-method and the post-method, we moreover compare the indicators $\eps_D^\nu$ and $\eps_1^\nu$ defined by:
\begin{equation}\label{eq:errUH1}
\begin{array}{rcl}
\eps_1^\nu&:=&\|\Grad(\Ical_{h,cg}\uvec-\uvec_h)\|_{\bbL^2(\Om)}/\|(\uvec,p)\|_{\Xcal,\nu},\\
\eps_D^\nu&:=&\|\div\uvec_h\|_{L^2(\Om)}/\|(\uvec,p)\|_{\Xcal,\nu}.
\end{array}
\end{equation} 
}
We first plot the results as error curves for both the velocity and the pressure as a function of the meshsize, and we give the numerically observed average convergence rates. Second, we report the elapsed CPU times for the {CR-method} and the {Post-method}. In the Tables, column "overhead" indicates the ratio between the cost of the second (post processing) step for the Post-method and the cost of the first step (CR-method). Finally, we plot the errors  of the two methods as a function of the elapsed CPU times of the resolution algorithm.

\textcolor{cyan}{For the Post-method "Post" plots represent the errors after a single iteration of the second step and "Post-8" plots represent the errors after iterating the second step eight times, updating ${p_h^{old}}$ at each new iteration, cf. Algorithm \eqref{eq:StokesFV2h-Iterative}.
\begin{rmk}
The use of discrete mixed formulations with $\bH(\div)$-conforming projection of the test function on the right-hand side leads to a pressure robust discrete velocity \cite{Link14,BLMS15}. It has been proven to be an efficient numerical solution to solve Stokes problem when the parameter $\nu$ is small, however the $\bH(\div)$-conforming discrete spaces must be chosen with care. 
The {Post}-method can also be used in that case to obtain an $\bH^1$-conforming reconstruction of the discrete velocity, see \S\ref{sec:RT}. 
\end{rmk}
}

\subsection{Regular manufactured solutions}\label{sec:RegSol}
We postulate that for the {EP-method} and the {Post-method}, $\eps_0^\nu(\uvec_h)\leq C_0\,h^2$, where $C_0$ is independent of the meshsize. The error estimates for the {CR-method} are given by \cite[Theorems 3, 4, 6]{CrRa73}.
Notice that the norm $\|(\cdot,\cdot)\|_{\Xcal,\nu}$ depends on the viscosity $\nu$ \eqref{eq:normX} and that the error estimates on the pressure and the velocity are linked.  We check the convergence rates and the viscosity dependency of the error estimates. The first test enters the framework of \S\ref{ss:cv-z-known}, see the estimate \eqref{eq:pressure-robust-estimate}.
\begin{itemize}
\item {\bf Test case with a linear velocity: Figures \ref{fig:XYFlow0}-\ref{fig:XYFlowCPU}, Tables \ref{tab:XYFlow}-\ref{tab:XYFlow-CPU}.}
\end{itemize}
We consider \textcolor{cyan}{Problem \eqref{eq:Stokes}$_{NH}$} with $\fvec=\grad p$, where: $\uvec=(-y,x)^T$ and $p=x^3+y^3-1/2$.
\\
The number of unknowns are reported on Table~\ref{tab:NoU} (left).\\
As expected in \eqref{eq:pressure-robust-estimate}, the {EP-method} returns $\eps_0^{\nu}(\uvec_h)=\Ocal(10^{-15})$, $\eps_0^\nu(p_h)=\Ocal(10^{-14})$ for $\nu=1$; and $\eps_0^{\nu}(\uvec_h)=\Ocal(10^{-18})$, $\eps_0^\nu(p_h)=\Ocal(10^{-17})$ for $\nu=10^{-6}$, so the errors are not reported. 

Figure \ref{fig:XYFlow0} (resp. \ref{fig:XYFlow}) shows the discrete error values $\eps_0^\nu(\uvec_h)$ (left) and \textcolor{cyan}{$\eps_0^\nu(p_h)$} (right) plotted against the meshsize for $\nu=1$ (resp. $\nu=10^{-6}$). \textcolor{cyan}{Notice that the Post-method improves the approximation initially computed with the CR-method.} 

\textcolor{cyan}{Figure \ref{fig:XYFlowH1Div} shows the discrete error values $\eps_1^\nu$ and $\eps_D^\nu$ plotted against the meshsize for $\nu=10^{-6}$. Notice that the ratio $\eps_D^\nu/\eps_1^\nu$ is of order $0.5$ for the Post-method with a single iteration and of order $0.1$ for the Post-method with eight iterations.}

Table \ref{tab:XYFlow} shows the average convergence rates for the velocity, $\tau_{\uvec}$ and the pressure, $\tau_p$. For the {CR-method}, the average convergence rate for the pressure, $\tau_p$ is better than expected, possibly because the source term is a polynomial of degree $3$ which gradient is numerically exactly integrated. 
 
Table \ref{tab:XYFlow-CPU} shows the CPU times for the CR-method and the Post-method with either one, or eight, iterations. With our stopping criterion, around $30$ iterations of the preconditioned conjugate gradient algorithm are needed to solve the pressure solver of algorithm \eqref{eq:Uzawa} with the {CR-method}: this is consistent with Remark~\ref{rem:cond}. By design, the Post-method which includes the CR-method as its first step requires more CPU time. However, we notice that the CPU time of the second (post processing) step, the overhead, is only a small fraction of the first step, and also that it decreases dramatically as the meshsize decreases. This can be explained by the fact that on the one hand, there are fewer unknowns and, on the other hand, algorithm \eqref{eq:Uzawa-Uh} is faster than algorithm \eqref{eq:Uzawa} (cf. Remark \ref{rem:cond}). 
For the {Post-method}, it seems worth doing eight iterations, especially when the meshsize is small. 

Figure \ref{fig:XYFlowCPU} shows the discrete error values $\eps_0^\nu(\uvec_h)$ (left) and $\eps_0^\nu(p_h)$ (right) plotted against the CPU time for $\nu=10^{-6}$. To reach $\eps_0^\nu(\uvec_h)\lesssim2\times10^{-5}$, the required CPU time is $0.5\,s$ for the Post-method with one iteration and less than $0.01\,s$ for the Post-method with eight iterations, to be compared with $20\,s$ for the CR-method.
\input{XYFlowFig0.tex}
\input{XYFlowFig.tex}
\input{XYFlowH1DivFig.tex}
\begin{table}[h]
\begin{center}
\begin{tabular}{c|c|c|c|c}
$\nu$&$\tau$& CR & Post & Post-8\\
\hline
\multirow{2}{*}{$1$}&$\tau_{\uvec}$&
$2.03$&$1.99$&$1.93$\\
&$\tau_p$&
$1.55$&$1.59$&$1.71$\\
\hline
\multirow{2}{*}{$10^{-6}$}&$\tau_{\uvec}$&
$2.03$&$1.97$&$2.53$\\
&$\tau_p$&
$1.55$&$1.68$& $1.78$
\end{tabular}
\caption{Linear velocity. Average convergence rates.}
\label{tab:XYFlow}
\end{center}
\end{table}
\begin{table}[h]
\begin{center}
\begin{tabular}{c|c|c|c|c|c}
$h$        & CPU CR     &  CPU Post & overhead &CPU Post-8 & overhead \\
\hline
$1.00\times 10^{-1}$& $4.34\times 10^{-3}$& $5.08\times 10^{-3}$& $17\,\%$ & $7.40\times 10^{-3}$ & $70\,\%$\\
$5.00\times 10^{-2}$& $2.32\times 10^{-2}$& $2.61\times 10^{-2}$& $13\,\%$ & $4.46\times 10^{-2}$ & $92\,\%$\\
$2.50\times 10^{-2}$& $2.51\times 10^{-1}$& $2.69\times 10^{-1}$& $7.2\,\%$ & $3.65\times 10^{-1}$ & $45\,\%$\\
$1.25\times 10^{-2}$& $3.00\times 10^{+0}$& $3.07\times 10^{+0}$& $2.3\,\%$ & $3.54\times 10^{+0}$ & $10\,\%$\\
$6.25\times 10^{-3}$& $4.89\times 10^{+1}$& $4.94\times 10^{+1}$& $1.0\,\%$  & $5.20\times 10^{+1}$ & $6.3\,\%$\\
\end{tabular}
\caption{Linear velocity, $\nu=10^{-6}$. CPU time (s).}
\label{tab:XYFlow-CPU}
\end{center}
\end{table}
\input{XYFlowCPUFig.tex}
\begin{itemize}
\item {\bf Test case  with a sinusoidal solution: Figures \ref{fig:SinFlow}-\ref{fig:SinFlowCPU6}, Tables \ref{tab:SinFlow-CPU}-\ref{tab:SinFlow}}
\end{itemize}
We consider \textcolor{cyan}{Problem \eqref{eq:Stokes}$_H$} with $\fvec=-\nu\Deltavec\uvec+\grad p$, where: \[
\uvec=\left(\begin{matrix}
(1-\cos(2\,\pi\,x))\,\sin(2\,\pi\,y)\\
(\cos(2\,\pi\,y)-1)\,\sin(2\,\pi\,x)
\end{matrix}\right)\mbox{ and }p=\sin(2\,\pi\,x)\,\sin(2\,\pi\,y).
\]
The number of unknowns are reported on Table~\ref{tab:NoU} (left).

Figure \ref{fig:SinFlow} shows $\eps_0^\nu (\uvec_h)$ (left) and $\eps_0^\nu (p_h)$ (right) plotted against the meshsize, for $\nu=1$. In that case, for the {Post-method}, there is no need doing eight iterations, so we give numerical results for computations with a single iteration only. 
The {Post-method} and the {EP-method} give similar velocity errors. The {CR-method} and the {Post-method} give similar pressure errors. For a given meshsize $h$, the velocity error is smaller for the {EP-method} and the {Post-method} than for the {CR-method}; and the pressure error is smaller for the {CR-method} than for the {EP-method}, but it is obtained at a higher cost. Notice that the {Post-method} reduces the velocity error without worsening the pressure error.

Figure \ref{fig:SinFlow6} shows $\eps_0^\nu (\uvec_h)$ (left) and $ \eps_0^\nu (p_h)$ (right) plotted against the meshsize, now for $\nu=10^{-6}$. Both post-processings improve the initial computation. Even with eight iterations, the overhead cost remains affordable, especially when the meshsize is small.
\\
The {EP-method} gives much smaller velocity and pressure errors than the {CR-method} since it is a pressure robust method. For a given meshsize $h$, we note that the {Post-method} allows to reduce the velocity error by a factor larger than $10$. The pressure error of the {Post-method} using one iteration is close to that of the {CR-method}, while it is greatly improved using eight iterations. 

\textcolor{cyan}{Figure \ref{fig:SinusFlowH1Div} shows the discrete error values $\eps_1^\nu$ and $\eps_D^\nu$ plotted against the meshsize for $\nu=10^{-6}$. Iterating several times yields again a significant decrease of the indicators.}

Table \ref{tab:SinFlow} shows the average convergence rates for the velocity, $\tau_{\uvec}$ and the pressure, $\tau_p$. In the case $\nu=1$, the average convergence rates are as expected. In the case $\nu=10^{-6}$, the average convergence rate for the velocity, $\tau_{\uvec}$, is better than expected for the {EP-method} and the {Post-method}, and the average convergence rate for the pressure, $\tau_p$, is better than expected for the {Post-method}, probably because the asymptotic convergence regime has not been reached.

Table \ref{tab:SinFlow-CPU} below shows the CPU times for the CR-method and the Post-method with either one, or eight, iterations. As we have used the same meshes, the computation times are similar to those in Table \ref{tab:XYFlow-CPU}, showing once again that the overhead cost decreases dramatically as the mesh size decreases.

Figure \ref{fig:SinFlowCPU6} shows $\eps_0^{\nu}(\uvec_h)$ (right) and $\eps_0^{\nu}(p_h)$ (left) plotted against the CPU time for the CR-method and the Post-method with a single iteration ("Post" plot) or eight iterations ("Post-8" plot). In the same way as for the first test, we note that to reach $\eps_0^{\nu}(\uvec_h)\lesssim 5\times 10^{-6}$, the required CPU time is $0.2\,s$ for the Post-method with one iteration and less than $0.02\,s$ for the Post-method with eight iterations, to be compared with the CR-method which does not reach this threshold.
\input{SinFlow0Fig.tex}
\input{SinFlowFig.tex}
\input{SinFlowH1DivFig.tex}
\begin{table}[h]
\begin{center}
\begin{tabular}{c|c|c|c|c|c}
$\nu$&$\tau$& CR & EP & Post & Post-8\\
\hline
\multirow{2}{*}{$1$}&$\tau_{\uvec}$&
$1.99$&$2.07$&$2.06$&$-$\\
&$\tau_p$&
$1.03$&$1.06$&$1.11$&$-$\\
\hline
\multirow{2}{*}{$10^{-6}$}&$\tau_{\uvec}$&
$2.05$&$2.17$&$2.47$&$2.36$\\
&$\tau_p$&
$1.31$&$1.07$&$1.47$& $1.74$
\end{tabular}
\caption{Sinusoidal velocity. Average convergence rates.}
\label{tab:SinFlow}
\end{center}
\end{table}
\begin{table}[h]
\begin{center}
\begin{tabular}{c|c|c|c|c|c}
$h$    & CPU CR & CPU Post & overhead & CPU Post-8 & overhead \\
\hline
$1.00\times 10^{-1}$& $3.85\times 10^{-3}$& $4.64\times 10^{-3}$& $21\,\%$& $7.06\times 10^{-3}$&$83\,\%$\\
$5.00\times 10^{-2}$& $2.22\times 10^{-2}$& $2.51\times 10^{-2}$& $13\,\%$& $4.27\times 10^{-2}$&$92\,\%$\\
$2.50\times 10^{-2}$& $2.43\times 10^{-1}$& $2.58\times 10^{-1}$& $6.2\,\%$& $3.45\times 10^{-1}$&$42\,\%$ \\
$1.25\times 10^{-2}$& $4.27\times 10^{+0}$& $4.34\times 10^{+0}$& $1.6\,\%$& $4.78\times 10^{+0}$&$12\,\%$\\
$6.25\times 10^{-3}$& $4.53\times 10^{+1}$& $4.57\times 10^{+1}$& $0.7\,\%$& $4.82\times 10^{+1}$&$6.4\,\%$\\
\end{tabular}
\caption{Sinusoidal velocity, $\nu=10^{-6}$. CPU time (s).}
\label{tab:SinFlow-CPU}
\end{center}
\end{table}
\input{SinFlowCPUFig}
\textcolor{cyan}{
\begin{itemize}
\item {\bf Comment on the algorithm}
\end{itemize}
At this stage, we noticed that eight iterations was a good compromise. In order to design an optimized algorithm, another stopping criterion could be set by comparing two successive computations of ${p_h^{old}}$, i.e. comparing $|\delta\ul{P}^n|_{\bbM}$ to $|\delta\ul{P}^{n-1}|_{\bbM}$. 
The choice of $\lambda$ could also be further optimized, depending on $\nu$ and on the meshsize.
}
\begin{itemize}
\item {\bf Some observations}
\end{itemize}
The coercive $\bP^1\times P^0$ formulation (EP-method) gives {\em pressure robust} results, the obvious limitation being that it requires to know explicitly the potential of the gradient part of the source term. If it is not known, the two step method greatly reduces the velocity error, compared with the calculation carried out using the Crouzeix-Raviart $\bP^1_{nc}\times P^0$ formulation (CR-method). Moreover, using the second (post processing) step {\em iteratively} greatly improves the initial result. Finally, the reduction factor is greater the smaller $\nu$ is.

\subsection{Low regularity manufactured solution: Figures \ref{fig:LowReg6h}-\ref{fig:LowRegCPU6h}, Tables \ref{tab:LowReg}-\ref{tab:LowReg-CPU}}\label{sec:LowRegSol}
Last, \textcolor{cyan}{we consider Problem \eqref{eq:Stokes}$_{NH}$} with a low regularity solution. Let $(\rho,\theta)$ be the polar coordinates centered in $(0.5,0.5)$. Let $\alpha=0.45$. We set $\fvec=-\nu\Deltavec{\uvec}+\grad p$ where $(\uvec,p)=(\rho^\alpha\evec_{\theta},\rho-\int_\Om\rho)$. Results are given for $\nu=10^{-6}$. 
The number of unknowns are reported on Table~\ref{tab:NoU} (right). We used a refined mesh around $(0.5,0.5)$, where the solution is of low regularity. 

Figure \ref{fig:LowReg6h} shows the discrete error values \textcolor{cyan}{$\eps_0^\nu (\uvec_h)$} (left) and \textcolor{cyan}{$ \eps_0^\nu (p_h)$} (right) plotted against the meshsize. We remark that {EP-method} shows far better results than the {CR-method}, and that the {Post-method} allows again to improve the approximation of the {CR-method}. 

\textcolor{cyan}{Figure \ref{fig:LowRegH1Div} shows the discrete error values $\eps_1^\nu$ (left) and $\eps_D^\nu$ (right) plotted against the meshsize.} We have $\eps_D^\nu/\eps_1^\nu\approx 0.5$ for the Post-method with a single iteration while $\eps_D^\nu/\eps_1^\nu\approx 0.1$ or the Post-method with eight iterations. 

Table \ref{tab:LowReg} shows the averaged convergence rates between the successive meshes. For the {EP-method} and the {CR-method}, we postulate that, asymptotically, $\tau_{\uvec}=1+\alpha$ and $\tau_p=\alpha$. In both cases, the obtained convergence rates are better than expected, which  suggest that the asymptotic convergence regime is not reached.

Table \ref{tab:LowReg-CPU} below shows the CPU times for the CR-method and the Post-method with eight iterations. Again, the overhead cost decreases sharply with the meshsize.

Figures \ref{fig:LowRegCPU6h} show the discrete error values $ \eps_0^\nu (\uvec_h)$ (left) and $ \eps_0^\nu (p_h)$ (right)  against the CPU time for the CR-method and the Post-method with eight iterations. To reach $\eps_0^\nu(\uvec_h)\lesssim 10^{-5}$, the required CPU time is $0.02\,s$ for the Post-method with eight iterations, to be compared with more than $100\,s$ for the CR-method.
\input{LowRegFig.tex}
\input{LowRegH1DivFig.tex}
\begin{table}[h]
\begin{center}\begin{tabular}{c|c|c|c}
$\tau$& CR & EP & Post-8\\
\hline
$\tau_{\uvec}$	& $2.0$	& $2.2$	& $2.8$\\
\hline
$\tau_p$		& $1.6$ & $0.7$ & $1.9$\\
\end{tabular}
\caption{low regularity velocity, $\nu=10^{-6}$. Averaged convergence rates.}
\label{tab:LowReg}
\end{center}
\end{table}
\\
\begin{table}[h]
\begin{center}
\begin{tabular}{c|c|c|c}
$h$ &CPU CR    & CPU Post-8 & overhead \\
\hline
$1.00\times 10^{-1}$& $1.12\times 10^{-2}$&	$1.92\times 10^{-2}$&$71\,\%$\\
$5.00\times 10^{-2}$& $2.77\times 10^{-1}$&	$3.31\times 10^{-1}$&$19\,\%$\\
$2.50\times 10^{-2}$& $3.15\times 10^{+0}$&	$3.40\times 10^{+0}$&$7.9\,\%$\\
$1.25\times 10^{-2}$& $1.06\times 10^{+1}$&	$1.20\times 10^{+1}$&$13\,\%$\\
$6.25\times 10^{-3}$& $1.32\times 10^{+2}$&	$1.39\times 10^{+2}$&$5\,\%$\\
\end{tabular}
\caption{Low regularity velocity, $\nu=10^{-6}$. CPU time (s).}
\label{tab:LowReg-CPU}
\end{center}
\end{table} 
\input{LowRegCPUFig.tex}
\subsection{\textcolor{cyan}{Numerical results using Raviart-Thomas projection}}\label{sec:RT}
The use of Crouzeix-Raviart $\bP^1_{nc}\times P^0$ formulation with $\bH(\div)$-conforming projection of the test function on the right-hand side leads to a pressure robust discrete velocity \cite{Link14,BLMS15}. In this case too, we can use the Post-method to obtain a precise $\bH^1$-conforming approximation of the discrete velocity. On Figures \ref{fig:XYFlowRT}-\ref{fig:LowRegRT6h}, we represent $\eps_0^\nu(\uvec_h)$ and $\eps_0^\nu(p_h)$ against the meshsize for our three tests. For the CR-method, we use the lowest order Raviart-Thomas projection \cite[\S3.2]{Link14} to compute ${p_h^{old}}$, which is then used in the Post-method with $\lambda=1$. For the first test (linear velocity), we obtain $\eps_0^\nu(\uvec_h)\lesssim 10^{-13}$ and $\eps_0^\nu(p_h)\lesssim 10^{-12}$, and regardless, we observe that the errors (close to machine precision for the CR-method) do not deteriorate when the Post-method is applied. The pressure error $\eps_0^\nu(p_h)$ remains unchanged for the other two tests. The approximation of the velocity is improved for all three tests.
\input{XYFlowFigRT.tex}
\input{SinFlowFigRT.tex}
\input{LowRegFigRT.tex}
\section{Conclusion}\label{sec-conclusion}
We proposed and analysed a new variational formulation of the Stokes problem based on  $T$-coercivity theory. This variational formulation is coercive, and can be discretized with the $\bP^k\times P_{disc}^{k-1}$ finite element for all $k\ge1$. To solve the linear system resulting from the discretization, we need to know the pressure, or at least some approximation of it, which in our numerical tests is the discrete pressure obtained using the classical non-conforming method with the Crouzeix-Raviart $\bP^1_{nc}\times P^0$ finite element. This two step method improves the numerical results by notably reducing the errors obtained after the use of the classical method, especially when the viscosity is small. More significantly, the two step method consistently outperforms the classical method in terms of precision with respect to CPU time.

%% file: XYFlowFig0.tex
\begin{figure}[h!]
\centering
\begin{tikzpicture}[scale=0.6]
\begin{axis}[
 height = 8.2cm,
 width = 8.2cm,
 xlabel = {$h$},
 ylabel = {$ \eps_0^\nu (\uvec_h)$},
 grid=both,
 major grid style={black!50},
 xmode=log, ymode=log,
 xmin=5e-3, xmax=1e-1,
 ymin=1e-9, ymax=1e-2,
 yticklabel style={
  /pgf/number format/fixed,
  /pgf/number format/precision=0},
 scaled y ticks=false,
 xticklabel style={
  /pgf/number format/fixed,
  /pgf/number format/precision=0},
 scaled y ticks=false,  
 legend style={at={(1,0)},anchor=south east},
 legend columns=2
],
\addplot table [x=h, y=eU0NC] {Resu/StokesXY-P1P0-nu0.txt};
\addlegendentry{CR}
\addplot table [x=h, y=eU0TS] {Resu/StokesXY-P1P0-nu0.txt};
\addlegendentry{Post}
\addplot table [x=h, y=eU0NS] {Resu/StokesXY-P1P0-nu0.txt};
\addlegendentry{Post-8}
\end{axis}       
\end{tikzpicture}
\begin{tikzpicture}[scale=0.6]
 \begin{axis}[
 height = 8.2cm, 
 width = 8.2cm,  
 xlabel = {$h$},
 ylabel = {$\eps_0^\nu(p_h)$},
 grid=both,
 major grid style={black!50},
 xmode=log, ymode=log,
 xmin=5e-3, xmax=1e-1,
 ymin=1e-5, ymax=1e-1,
 yticklabel style={
     /pgf/number format/fixed,
     /pgf/number format/precision=0
 },
 scaled y ticks=false,
 xticklabel style={
     /pgf/number format/fixed,
     /pgf/number format/precision=0
 },
 scaled y ticks=false,  
 legend style={at={(1,0)},anchor=south east},
 legend columns=2
 ],
 \addplot table [x=h, y=eP0NC] {Resu/StokesXY-P1P0-nu0.txt};
 \addlegendentry{CR}
 \addplot table [x=h, y=eP0TS] {Resu/StokesXY-P1P0-nu0.txt};
 \addlegendentry{Post}
 \addplot table [x=h, y=eP0NS] {Resu/StokesXY-P1P0-nu0.txt};
 \addlegendentry{Post-8}
 \end{axis}      
\end{tikzpicture}
\caption{Linear velocity, $\nu=1$.
Plots of $ \eps_0^\nu (\uvec_h)$ (left) and $\eps_0^\nu(p_h)$ (right).}
\label{fig:XYFlow0}
\end{figure}

%% file: XYFlowFig.tex
\begin{figure}[h!]
\centering
\begin{tikzpicture}[scale=0.6]
\begin{axis}[
 height = 8.2cm,
 width = 8.2cm,
 xlabel = {$h$},
 ylabel = {$ \eps_0^\nu (\uvec_h)$},
 grid=both,
 major grid style={black!50},
 xmode=log, ymode=log,
 xmin=5e-3, xmax=1e-1,
 ymin=1e-9, ymax=1e-2,
 yticklabel style={
  /pgf/number format/fixed,
  /pgf/number format/precision=0},
 scaled y ticks=false,
 xticklabel style={
  /pgf/number format/fixed,
  /pgf/number format/precision=0},
 scaled y ticks=false,  
 legend style={at={(1,0)},anchor=south east},
 legend columns=2
],
\addplot table [x=h, y=eU0NC] {Resu/StokesXY-P1P0-nu6.txt};
\addlegendentry{CR}
\addplot table [x=h, y=eU0TS] {Resu/StokesXY-P1P0-nu6.txt};
\addlegendentry{Post}
\addplot table [x=h, y=eU0NS] {Resu/StokesXY-P1P0-nu6.txt};
\addlegendentry{Post-8}
\end{axis}       
\end{tikzpicture}
\begin{tikzpicture}[scale=0.6]
 \begin{axis}[
 height = 8.2cm, 
 width = 8.2cm,  
 xlabel = {$h$},
 ylabel = {$\eps_0^\nu(p_h)$},
 grid=both,
 major grid style={black!50},
 xmode=log, ymode=log,
 xmin=5e-3, xmax=1e-1,
 ymin=1e-5, ymax=1e-1,
 yticklabel style={
     /pgf/number format/fixed,
     /pgf/number format/precision=0
 },
 scaled y ticks=false,
 xticklabel style={
     /pgf/number format/fixed,
     /pgf/number format/precision=0
 },
 scaled y ticks=false,  
 legend style={at={(1,0)},anchor=south east},
 legend columns=2
 ],
 \addplot table [x=h, y=eP0NC] {Resu/StokesXY-P1P0-nu6.txt};
 \addlegendentry{CR}
 \addplot table [x=h, y=eP0TS] {Resu/StokesXY-P1P0-nu6.txt};
 \addlegendentry{Post}
 \addplot table [x=h, y=eP0NS] {Resu/StokesXY-P1P0-nu6.txt};
 \addlegendentry{Post-8}
 \end{axis}      
\end{tikzpicture}
\caption{Linear velocity, $\nu=10^{-6}$.
Plots of $ \eps_0^\nu (\uvec_h)$ (left) and $\eps_0^\nu(p_h)$ (right).}
\label{fig:XYFlow}
\end{figure}

%% file: XYFlowH1DivFig.tex
\begin{figure}[h!]
\centering
\begin{tikzpicture}[scale=0.6]
\begin{axis}[
 height = 8.2cm,
 width = 8.2cm,
 xlabel = {$h$},
 ylabel = {$\eps_1^\nu$},
 grid=both,
 major grid style={black!50},
 xmode=log, ymode=log,
 xmin=5e-3, xmax=1e-1,
 ymin=1e-8, ymax=1e-2,
 yticklabel style={
  /pgf/number format/fixed,
  /pgf/number format/precision=0},
 scaled y ticks=false,
 xticklabel style={
  /pgf/number format/fixed,
  /pgf/number format/precision=0},
 scaled y ticks=false,  
 legend style={at={(1,0)},anchor=south east},
 legend columns=2
],
\addplot table [x=h, y=eU1TS] {Resu/StokesXY-P1P0-nu6.txt};
\addlegendentry{Post, $\eps_1^\nu$}
\addplot table [x=h, y=eU1NS] {Resu/StokesXY-P1P0-nu6.txt};
\addlegendentry{Post-8, $\eps_1^\nu$}
\end{axis}       
\end{tikzpicture}
\begin{tikzpicture}[scale=0.6]
 \begin{axis}[
 height = 8.2cm, 
 width = 8.2cm,  
 xlabel = {$h$},
 ylabel = {$\eps_D^\nu$},
 grid=both,
 major grid style={black!50},
 xmode=log, ymode=log,
 xmin=5e-3, xmax=1e-1,
 ymin=1e-8, ymax=1e-2,
 yticklabel style={
     /pgf/number format/fixed,
     /pgf/number format/precision=0
 },
 scaled y ticks=false,
 xticklabel style={
     /pgf/number format/fixed,
     /pgf/number format/precision=0
 },
 scaled y ticks=false,  
 legend style={at={(1,0)},anchor=south east},
 legend columns=2
 ],
\addplot table [x=h, y=eUDTS] {Resu/StokesXY-P1P0-nu6.txt};
\addlegendentry{Post, $\eps_D^\nu$}
\addplot table [x=h, y=eUDNS] {Resu/StokesXY-P1P0-nu6.txt};
\addlegendentry{Post-8, $\eps_D^\nu$}
 \end{axis}      
\end{tikzpicture}
\caption{Linear velocity, $\nu=10^{-6}$.
Plots of $\eps_1^\nu$  (left) and $\eps_D^\nu$ (right).}
\label{fig:XYFlowH1Div}
\end{figure}
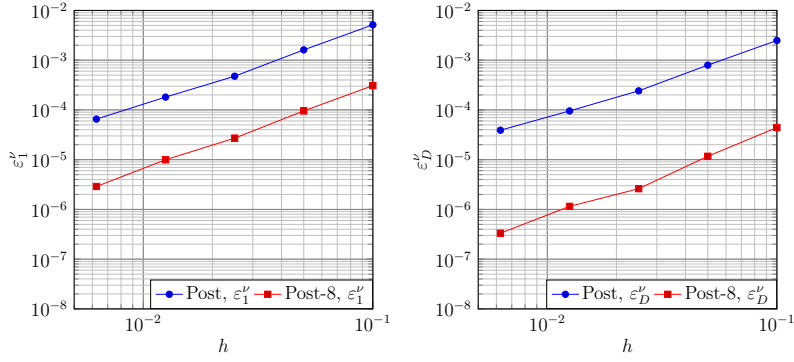

%% file: XYFlowCPUFig.tex
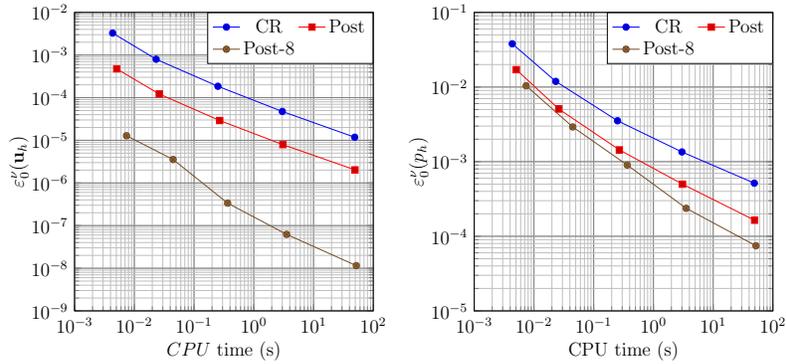
\begin{figure}[h!]
\centering
\begin{tikzpicture}[scale=0.6]
\begin{axis}[
 height = 8.2cm,
 width = 8.2cm,
 xlabel = {$CPU$ time (s)},
 ylabel = {$ \eps_0^\nu (\uvec_h)$},
 grid=both,
 major grid style={black!50},
 xmode=log, ymode=log,
 xmin=1e-3, xmax=1e+2,
 ymin=1e-9, ymax=1e-2,
 yticklabel style={
  /pgf/number format/fixed,
  /pgf/number format/precision=0},
 scaled y ticks=false,
 xticklabel style={
  /pgf/number format/fixed,
  /pgf/number format/precision=0},
 scaled y ticks=false,  
 legend style={at={(1,1)}},
 legend columns=2
],
\addplot table [x=tpsNC, y=eU0NC] {Resu/StokesXY-P1P0-nu6.txt};
\addlegendentry{CR}
\addplot table [x=tpsTS, y=eU0TS] {Resu/StokesXY-P1P0-nu6.txt};
\addlegendentry{Post}
\addplot table [x=tpsNS, y=eU0NS] {Resu/StokesXY-P1P0-nu6.txt};
\addlegendentry{Post-8}
\end{axis}       
\end{tikzpicture}
\begin{tikzpicture}[scale=0.6]
 \begin{axis}[
 height = 8.2cm,
 width = 8.2cm,
 ylabel = {$\eps_0^\nu(p_h)$},
 xlabel = {CPU time (s)},
 grid=both,
 major grid style={black!50},
 xmode=log, ymode=log,
 xmin=1e-3, xmax=1e+2,
 ymin=1e-5, ymax=1e-1,
 yticklabel style={
     /pgf/number format/fixed,
     /pgf/number format/precision=0
 },
 scaled y ticks=false,
 xticklabel style={
     /pgf/number format/fixed,
     /pgf/number format/precision=0
 },
 scaled y ticks=false,
 legend style={at={(1,1)}},
 legend columns=2
 ],
 \addplot table [x=tpsNC, y=eP0NC] {Resu/StokesXY-P1P0-nu6.txt};
 \addlegendentry{CR}
 \addplot table [x=tpsTS, y=eP0TS] {Resu/StokesXY-P1P0-nu6.txt};
 \addlegendentry{Post}
 \addplot table [x=tpsNS, y=eP0NS] {Resu/StokesXY-P1P0-nu6.txt};
 \addlegendentry{Post-8}
 \end{axis}      
\end{tikzpicture}
\caption{Linear velocity, $\nu=10^{-6}$.
Plots of $ \eps_0^\nu (\uvec_h)$ (left) and $\eps_0^\nu(p_h)$ (right) against CPU time (s).}
\label{fig:XYFlowCPU}
\end{figure}

%% file: SinFlow0Fig.tex
\begin{figure}[h!]
\centering
\begin{tikzpicture}[scale=0.6]
 \begin{axis}[
 height = 8.2cm,
 width = 8.2cm,
 xlabel = {$h$},
 ylabel = {$\eps_0^\nu(\uvec_h)$},
 grid=both,
 major grid style={black!50},
 xmode=log, ymode=log,
 xmin=5e-3, xmax=1e-1,
 ymin=1e-6, ymax=1e-2,
 yticklabel style={
     /pgf/number format/fixed,
     /pgf/number format/precision=0
 },
 scaled y ticks=false,
 xticklabel style={
     /pgf/number format/fixed,
     /pgf/number format/precision=0
 },
 scaled y ticks=false,  
 legend style={at={(1,0)},anchor=south east},
 legend columns=2
 ],
 \addplot table [x=h, y=eU0NC] {Resu/StokesSinus-P1P0-nu0.txt};
 \addlegendentry{CR}
 \addplot table [x=h, y=eU0EP] {Resu/StokesSinus-P1P0-nu0.txt};
 \addlegendentry{EP}
 \addplot table [x=h, y=eU0TS] {Resu/StokesSinus-P1P0-nu0.txt};
 \addlegendentry{Post}
 \end{axis}
\end{tikzpicture}
\begin{tikzpicture}[scale=0.6]
 \begin{axis}[
 height = 8.2cm,
 width = 8.2cm,
 xlabel = {$h$},
 ylabel = {$\eps_0^\nu(p_h)$},
 grid=both,
 major grid style={black!50},
 xmode=log, ymode=log,
 xmin=5e-3, xmax=1e-1,
 ymin=1e-3, ymax=1e-1,
 yticklabel style={
     /pgf/number format/fixed,
     /pgf/number format/precision=0
 },
 scaled y ticks=false,
 xticklabel style={
     /pgf/number format/fixed,
     /pgf/number format/precision=0
 },
 scaled y ticks=false,  
 legend style={at={(1,0)},anchor=south east},
 legend columns=2
 ],
 \addplot table [x=h, y=eP0NC] {Resu/StokesSinus-P1P0-nu0.txt};
 \addlegendentry{CR}
 \addplot table [x=h, y=eP0EP] {Resu/StokesSinus-P1P0-nu0.txt};
 \addlegendentry{EP}
 \addplot table [x=h, y=eP0TS] {Resu/StokesSinus-P1P0-nu0.txt};
 \addlegendentry{Post}
 \end{axis}      
\end{tikzpicture}
\caption{Sinusoidal velocity, $\nu=1$.
Plots of $\eps_0^\nu(\uvec_h)$ (left) and $\eps_0^\nu(p_h)$ (right).}
\label{fig:SinFlow}
\end{figure}
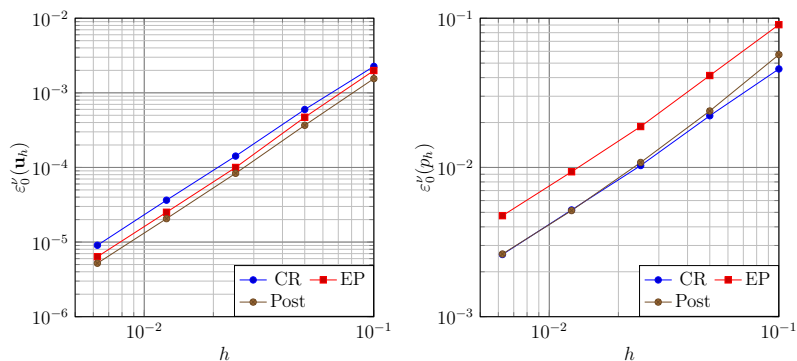

%% file: SinFlowFig.tex
\begin{figure}[h!]
\centering
\begin{tikzpicture}[scale=0.6]
 \begin{axis}[
 height = 8.2cm,
 width = 8.2cm,
 xlabel = {$h$},
 ylabel = {$\eps_0^\nu(\uvec_h)$},
 grid=both,
 major grid style={black!50},
 xmode=log, ymode=log,
 xmin=5e-3, xmax=1e-1,
 ymin=1e-11, ymax=1e-2,
 yticklabel style={
     /pgf/number format/fixed,
     /pgf/number format/precision=0
 },
 scaled y ticks=false,
 xticklabel style={
     /pgf/number format/fixed,
     /pgf/number format/precision=0
 },
 scaled y ticks=false,  
 legend style={at={(1,0)},anchor=south east},
 legend columns=2
 ],
 \addplot table [x=h, y=eU0NC] {Resu/StokesSinus-P1P0-nu6.txt};
 \addlegendentry{CR}
 \addplot table [x=h, y=eU0EP] {Resu/StokesSinus-P1P0-nu6.txt};
 \addlegendentry{EP}
 \addplot table [x=h, y=eU0TS] {Resu/StokesSinus-P1P0-nu6.txt};
 \addlegendentry{Post}
  \addplot table [x=h, y=eU0NS] {Resu/StokesSinus-P1P0-nu6.txt};
 \addlegendentry{Post-8}
 \end{axis}
\end{tikzpicture}
\begin{tikzpicture}[scale=0.6]
 \begin{axis}[
 height = 8.2cm,
 width = 8.2cm,
 xlabel = {$h$},
 ylabel = {$\eps_0^\nu (p_h)$},
 grid=both,
 major grid style={black!50},
 xmode=log, ymode=log,
 xmin=5e-3, xmax=1e-1,
 ymin=1e-8, ymax=1e-1,
 yticklabel style={
     /pgf/number format/fixed,
     /pgf/number format/precision=0
 },
 scaled y ticks=false,
 xticklabel style={
     /pgf/number format/fixed,
     /pgf/number format/precision=0
 },
 scaled y ticks=false,  
 legend style={at={(1,0)},anchor=south east},
 legend columns=2
 ],
 \addplot table [x=h, y=eP0NC] {Resu/StokesSinus-P1P0-nu6.txt};
 \addlegendentry{CR}
 \addplot table [x=h, y=eP0EP] {Resu/StokesSinus-P1P0-nu6.txt};
 \addlegendentry{EP}
 \addplot table [x=h, y=eP0TS] {Resu/StokesSinus-P1P0-nu6.txt};
 \addlegendentry{Post}
  \addplot table [x=h, y=eP0NS] {Resu/StokesSinus-P1P0-nu6.txt};
 \addlegendentry{Post-8}
 \end{axis}      
\end{tikzpicture}
\caption{Sinusoidal velocity, $\nu=10^{-6}$.
Plots of $\eps_0^\nu(\uvec_h)$ (left) and $\eps_0^\nu(p_h)$ (right).}
\label{fig:SinFlow6}
\end{figure}

%% file: SinFlowH1DivFig.tex
\begin{figure}[h!]
\centering
\begin{tikzpicture}[scale=0.6]
\begin{axis}[
 height = 8.2cm,
 width = 8.2cm,
 xlabel = {$h$},
 ylabel = {$\eps_1^\nu$},
 grid=both,
 major grid style={black!50},
 xmode=log, ymode=log,
 xmin=5e-3, xmax=1e-1,
 ymin=1e-7, ymax=1e-2,
 yticklabel style={
  /pgf/number format/fixed,
  /pgf/number format/precision=0},
 scaled y ticks=false,
 xticklabel style={
  /pgf/number format/fixed,
  /pgf/number format/precision=0},
 scaled y ticks=false,  
 legend style={at={(1,0)},anchor=south east},
 legend columns=2
],
\addplot table [x=h, y=eU1TS] {Resu/StokesSinus-P1P0-nu6.txt};
\addlegendentry{Post, $\eps_1^\nu$}
\addplot table [x=h, y=eU1NS] {Resu/StokesSinus-P1P0-nu6.txt};
\addlegendentry{Post-8, $\eps_1^\nu$}
\end{axis}       
\end{tikzpicture}
\begin{tikzpicture}[scale=0.6]
 \begin{axis}[
 height = 8.2cm, 
 width = 8.2cm,  
 xlabel = {$h$},
 ylabel = {$\eps_D^\nu$},
 grid=both,
 major grid style={black!50},
 xmode=log, ymode=log,
 xmin=5e-3, xmax=1e-1,
 ymin=1e-7, ymax=1e-2,
 yticklabel style={
     /pgf/number format/fixed,
     /pgf/number format/precision=0
 },
 scaled y ticks=false,
 xticklabel style={
     /pgf/number format/fixed,
     /pgf/number format/precision=0
 },
 scaled y ticks=false,  
 legend style={at={(1,0)},anchor=south east},
 legend columns=2
 ],
\addplot table [x=h, y=eUDTS] {Resu/StokesSinus-P1P0-nu6.txt};
\addlegendentry{Post, $\eps_D^\nu$}
\addplot table [x=h, y=eUDNS] {Resu/StokesSinus-P1P0-nu6.txt};
\addlegendentry{Post-8, $\eps_D^\nu$}
 \end{axis}      
\end{tikzpicture}
\caption{Sinusoidal velocity, $\nu=10^{-6}$.
Plots of $\eps_1^\nu$ (left) and $\eps_D^\nu$ (right).}
\label{fig:SinusFlowH1Div}
\end{figure}
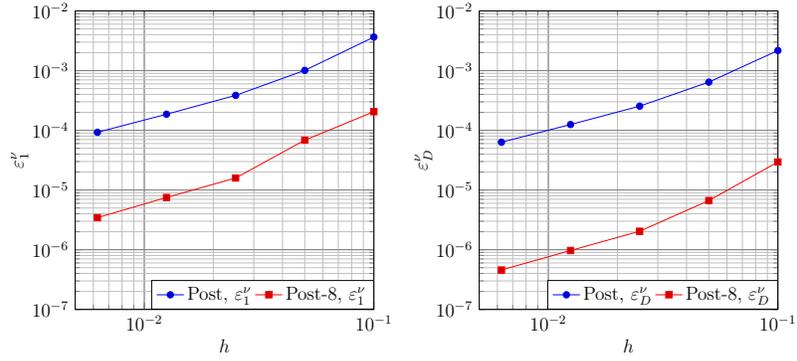

%% file: SinFlowCPUFig.tex
\begin{figure}[h!]
\centering
\begin{tikzpicture}[scale=0.6]
 \begin{axis}[
 height = 8.2cm, 
 width = 8.2cm,  
 ylabel = {$\eps_0^\nu(\uvec_h)$},
 xlabel = {CPU time (s)},
 grid=both,
 major grid style={black!50},
 xmode=log, ymode=log,
 ymin=1e-9, ymax=1e-2,
 xmin=1e-3, xmax=1e+2,
 yticklabel style={
     /pgf/number format/fixed,
     /pgf/number format/precision=0
 },
 scaled y ticks=false,
 xticklabel style={
     /pgf/number format/fixed,
     /pgf/number format/precision=0
 },
 scaled y ticks=false,  
 legend style={at={(0,0)},anchor=south west},
 legend columns=2
 ],
 \addplot table [x=tpsNC, y=eU0NC] {Resu/StokesSinus-P1P0-nu6.txt};
 \addlegendentry{CR}
 \addplot table [x=tpsTS, y=eU0TS] {Resu/StokesSinus-P1P0-nu6.txt};
 \addlegendentry{Post}
  \addplot table [x=tpsNS, y=eU0NS] {Resu/StokesSinus-P1P0-nu6.txt};
 \addlegendentry{Post-8}
 \end{axis}
\end{tikzpicture}
\begin{tikzpicture}[scale=0.6]
 \begin{axis}[
 height = 8.2cm, 
 width = 8.2cm,  
 ylabel = {$\eps_0^\nu(p_h)$},
 xlabel = {CPU time (s)},
 grid=both,
 major grid style={black!50},
 xmode=log, ymode=log,
 ymin=1e-5, ymax=1e-1,
 xmin=1e-3, xmax=1e+2,
 yticklabel style={
     /pgf/number format/fixed,
     /pgf/number format/precision=0
 },
 scaled y ticks=false,
 xticklabel style={
     /pgf/number format/fixed,
     /pgf/number format/precision=0
 },
 scaled y ticks=false,  
 legend style={at={(0,0)},anchor=south west},
 legend columns=2
 ],
 \addplot table [x=tpsNC, y=eP0NC] {Resu/StokesSinus-P1P0-nu6.txt};
 \addlegendentry{CR}
 \addplot table [x=tpsTS, y=eP0TS] {Resu/StokesSinus-P1P0-nu6.txt};
 \addlegendentry{Post}
  \addplot table [x=tpsNS, y=eP0NS] {Resu/StokesSinus-P1P0-nu6.txt};
 \addlegendentry{Post-8}
 \end{axis}      
\end{tikzpicture}
\caption{Sinusoidal velocity, $\nu=10^{-6}$.
Plots of $\eps_0^\nu(\uvec_h)$ (left) and $\eps_0^\nu(p_h)$ (right) against CPU time (s).}
\label{fig:SinFlowCPU6}
\end{figure}
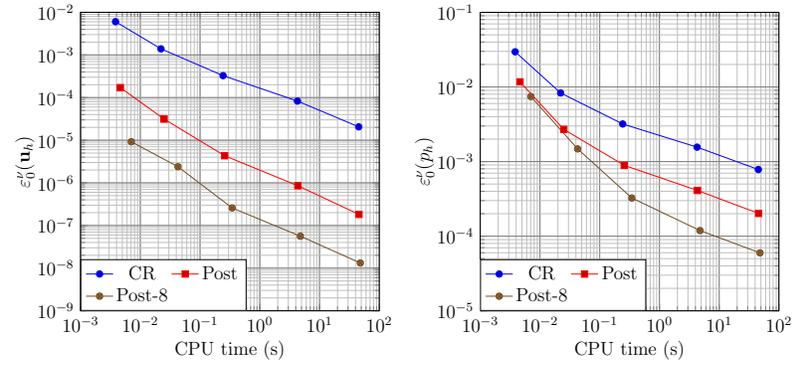

%% file: LowRegFig.tex
\begin{figure}[h!]
\centering
\begin{tikzpicture}[scale=0.6]
 \begin{axis}[
 height = 8.2cm,
 width = 8.2cm,
 xlabel = {$h$},
 ylabel = {$\eps_0^\nu(\uvec_h)$},
 grid=both,
 major grid style={black!50},
 xmode=log, ymode=log,
 xmin=5e-3, xmax=1e-1,
 ymin=1e-11, ymax=1e-2,
 yticklabel style={
     /pgf/number format/fixed,
     /pgf/number format/precision=0
 },
 scaled y ticks=false,
 xticklabel style={
     /pgf/number format/fixed,
     /pgf/number format/precision=0
 },
 scaled y ticks=false,  
 legend style={at={(1,0)},anchor=south east},
 legend columns=2
 ],
 \addplot table [x=h, y=eU0NC] {Resu/StokesVortex-P1P0-nu6.txt};
 \addlegendentry{CR}
 \addplot table [x=h, y=eU0EP] {Resu/StokesVortex-P1P0-nu6.txt};
 \addlegendentry{EP}
 \addplot table [x=h, y=eU0TS] {Resu/StokesVortex-P1P0-nu6.txt};
 \addlegendentry{Post}
 \addplot table [x=h, y=eU0NS] {Resu/StokesVortex-P1P0-nu6.txt};
 \addlegendentry{Post-8}
 \end{axis}
\end{tikzpicture}
\begin{tikzpicture}[scale=0.6]
 \begin{axis}[
 height = 8.2cm, 
 width = 8.2cm,  
 xlabel = {$h$},
 ylabel = {$\eps_0^\nu(p_h)$},
 grid=both,
 major grid style={black!50},
 xmode=log, ymode=log,
 xmin=5e-3, xmax=1e-1,
 ymin=1e-8, ymax=1e-1,
 yticklabel style={
     /pgf/number format/fixed,
     /pgf/number format/precision=0
 },
 scaled y ticks=false,
 xticklabel style={
     /pgf/number format/fixed,
     /pgf/number format/precision=0
 },
 scaled y ticks=false,  
 legend style={at={(1,0)},anchor=south east},
 legend columns=2
 ],
 \addplot table [x=h, y=eP0NC] {Resu/StokesVortex-P1P0-nu6.txt};
 \addlegendentry{CR}
 \addplot table [x=h, y=eP0EP] {Resu/StokesVortex-P1P0-nu6.txt};
 \addlegendentry{EP}
 \addplot table [x=h, y=eP0TS] {Resu/StokesVortex-P1P0-nu6.txt};
 \addlegendentry{Post}
 \addplot table [x=h, y=eP0NS] {Resu/StokesVortex-P1P0-nu6.txt};
 \addlegendentry{Post-8}
 \end{axis}      
\end{tikzpicture}
\caption{Low regularity velocity, $\nu=10^{-6}$.
Plots of $\eps_0^\nu(\uvec_h)$ (left) and $\eps_0^\nu (p_h)$ (right).}
\label{fig:LowReg6h}
\end{figure}

%% file: LowRegH1DivFig.tex
\begin{figure}[h!]
\centering
\begin{tikzpicture}[scale=0.6]
\begin{axis}[
 height = 8.2cm,
 width = 8.2cm,
 xlabel = {$h$},
 ylabel = {$\eps_1^\nu$},
 grid=both,
 major grid style={black!50},
 xmode=log, ymode=log,
 xmin=5e-3, xmax=1e-1,
 ymin=1e-7, ymax=1e-2,
 yticklabel style={
  /pgf/number format/fixed,
  /pgf/number format/precision=0},
 scaled y ticks=false,
 xticklabel style={
  /pgf/number format/fixed,
  /pgf/number format/precision=0},
 scaled y ticks=false,  
 legend style={at={(1,0)},anchor=south east},
 legend columns=2
],
\addplot table [x=h, y=eU1TS] {Resu/StokesVortex-P1P0-nu6.txt};
\addlegendentry{Post}
\addplot table [x=h, y=eU1NS] {Resu/StokesVortex-P1P0-nu6.txt};
\addlegendentry{Post-8}
\end{axis}       
\end{tikzpicture}
\begin{tikzpicture}[scale=0.6]
 \begin{axis}[
 height = 8.2cm, 
 width = 8.2cm,  
 xlabel = {$h$},
 ylabel = {$\eps_D^\nu$},
 grid=both,
 major grid style={black!50},
 xmode=log, ymode=log,
 xmin=5e-3, xmax=1e-1,
 ymin=1e-7, ymax=1e-2,
 yticklabel style={
     /pgf/number format/fixed,
     /pgf/number format/precision=0
 },
 scaled y ticks=false,
 xticklabel style={
     /pgf/number format/fixed,
     /pgf/number format/precision=0
 },
 scaled y ticks=false,  
 legend style={at={(1,0)},anchor=south east},
 legend columns=2
 ],
\addplot table [x=h, y=eUDTS] {Resu/StokesVortex-P1P0-nu6.txt};
\addlegendentry{Post}
\addplot table [x=h, y=eUDNS] {Resu/StokesVortex-P1P0-nu6.txt};
\addlegendentry{Post-8}
 \end{axis}      
\end{tikzpicture}
\caption{Low regularity velocity, $\nu=10^{-6}$.
Plots of $\eps_1^\nu$ (left) and $\eps_D^\nu$ (right).}
\label{fig:LowRegH1Div}
\end{figure}
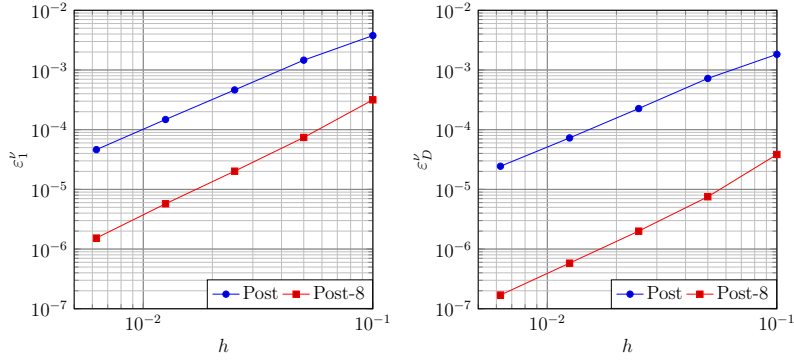

%% file: LowRegCPUFig.tex
\begin{figure}[h!]
\centering
\begin{tikzpicture}[scale=0.6]
 \begin{axis}[
 height = 8.2cm,
 width = 8.2cm,
 ylabel = {$\eps_0^\nu(\uvec_h)$},
 xlabel = {CPU time (s)},
 grid=both,
 major grid style={black!50},
 xmode=log, ymode=log,
 ymin=1e-9, ymax=1e-2,
 xmin=1e-2, xmax=1e+3,
 yticklabel style={
     /pgf/number format/fixed,
     /pgf/number format/precision=0
 },
 scaled y ticks=false,
 xticklabel style={
     /pgf/number format/fixed,
     /pgf/number format/precision=0
 },
 scaled y ticks=false,  
 legend style={at={(0,0)},anchor=south west},
 legend columns=2
 ],
 \addplot table [x=tpsNC, y=eU0NC] {Resu/StokesVortex-P1P0-nu6.txt};
 \addlegendentry{CR}
 \addplot table [x=tpsNS, y=eU0NS] {Resu/StokesVortex-P1P0-nu6.txt};
 \addlegendentry{Post-8}
 \end{axis}
\end{tikzpicture}
\begin{tikzpicture}[scale=0.6]
 \begin{axis}[
 height = 8.2cm,
 width = 8.2cm,
 ylabel = {$\eps_0^\nu(p_h)$},
 xlabel = {CPU time (s)},
 grid=both,
 major grid style={black!50},
 xmode=log, ymode=log,
 ymin=1e-5, ymax=1e-1,
 xmin=1e-2, xmax=1e+3,
 yticklabel style={
     /pgf/number format/fixed,
     /pgf/number format/precision=0
 },
 scaled y ticks=false,
 xticklabel style={
     /pgf/number format/fixed,
     /pgf/number format/precision=0
 },
 scaled y ticks=false,  
 legend style={at={(0,0)},anchor=south west},
 legend columns=2
 ],
 \addplot table [x=tpsNC, y=eP0NC] {Resu/StokesVortex-P1P0-nu6.txt};
 \addlegendentry{CR}
 \addplot table [x=tpsNS, y=eP0NS] {Resu/StokesVortex-P1P0-nu6.txt};
 \addlegendentry{Post-8}
 \end{axis}      
\end{tikzpicture}
\caption{Low regularity velocity, $\nu=10^{-6}$. Plots of $\eps_0^\nu(\uvec_h)$ (left) and $\eps_0^\nu(p_h)$ (right) against CPU time (s).}
\label{fig:LowRegCPU6h}
\end{figure}
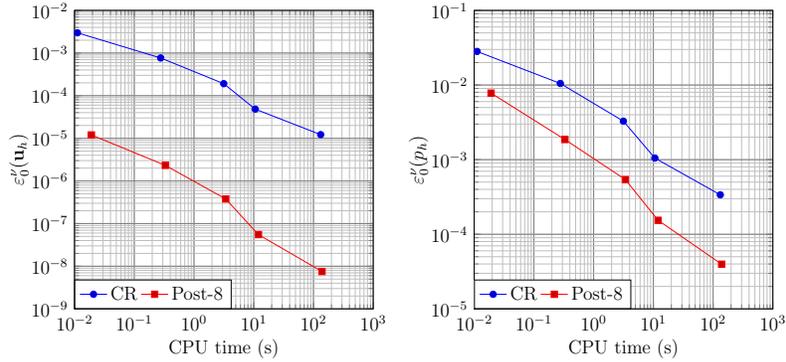

%% file: XYFlowFigRT.tex
\begin{figure}[h!]
\centering
\begin{tikzpicture}[scale=0.6]
\begin{axis}[
 height = 8.2cm,
 width = 8.2cm,
 xlabel = {$h$},
 ylabel = {$ \eps_0^\nu (\uvec_h)$},
 grid=both,
 major grid style={black!50},
 xmode=log, ymode=log,
 xmin=5e-3, xmax=1e-1,
 ymin=1e-15, ymax=1e-12,
 yticklabel style={
  /pgf/number format/fixed,
  /pgf/number format/precision=0},
 scaled y ticks=false,
 xticklabel style={
  /pgf/number format/fixed,
  /pgf/number format/precision=0},
 scaled y ticks=false,  
 legend style={at={(1,0)},anchor=south east},
 legend columns=2
],
\addplot table [x=h, y=eU0RT] {Resu/StokesXYRT-P1P0-nu6.txt};
\addlegendentry{CR}
\addplot table [x=h, y=eU0TS] {Resu/StokesXYRT-P1P0-nu6.txt};
\addlegendentry{Post}
\end{axis}       
\end{tikzpicture}
\begin{tikzpicture}[scale=0.6]
 \begin{axis}[
 height = 8.2cm, 
 width = 8.2cm,  
 xlabel = {$h$},
 ylabel = {$\eps_0^\nu(p_h)$},
 grid=both,
 major grid style={black!50},
 xmode=log, ymode=log,
 xmin=5e-3, xmax=1e-1,
 ymin=1e-15, ymax=1e-12,
 yticklabel style={
     /pgf/number format/fixed,
     /pgf/number format/precision=0
 },
 scaled y ticks=false,
 xticklabel style={
     /pgf/number format/fixed,
     /pgf/number format/precision=0
 },
 scaled y ticks=false,  
 legend style={at={(1,0)},anchor=south east},
 legend columns=2
 ],
 \addplot table [x=h, y=eU0RT] {Resu/StokesXYRT-P1P0-nu6.txt};
 \addlegendentry{CR}
 \addplot table [x=h, y=eP0TS] {Resu/StokesXYRT-P1P0-nu6.txt};
 \addlegendentry{Post}
 \end{axis}      
\end{tikzpicture}
\caption{Linear velocity, $\nu=10^{-6}$.
Plots of $ \eps_0^\nu (\uvec_h)$ (left) and $\eps_0^\nu(p_h)$ (right).}
\label{fig:XYFlowRT}
\end{figure}
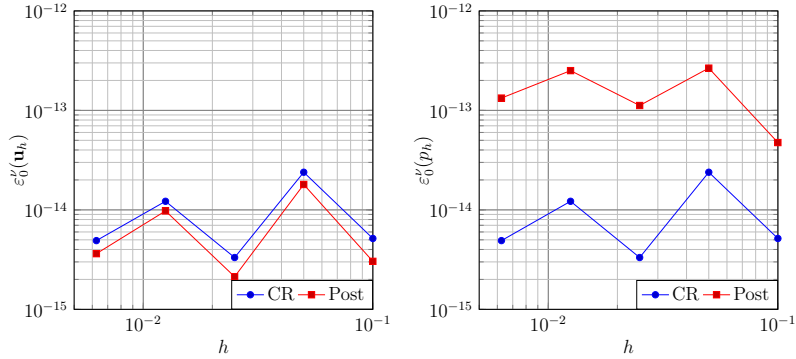

%% file: SinFlowFigRT.tex
\begin{figure}[h!]
\centering
\begin{tikzpicture}[scale=0.6]
 \begin{axis}[
 height = 8.2cm,
 width = 8.2cm,
 xlabel = {$h$},
 ylabel = {$\eps_0^\nu(\uvec_h)$},
 grid=both,
 major grid style={black!50},
 xmode=log, ymode=log,
 xmin=5e-3, xmax=1e-1,
 ymin=1e-11, ymax=1e-5,
 yticklabel style={
     /pgf/number format/fixed,
     /pgf/number format/precision=0
 },
 scaled y ticks=false,
 xticklabel style={
     /pgf/number format/fixed,
     /pgf/number format/precision=0
 },
 scaled y ticks=false,  
 legend style={at={(1,0)},anchor=south east},
 legend columns=2
 ],
 \addplot table [x=h, y=eU0RT] {Resu/StokesSinusRT-P1P0-nu6.txt};
 \addlegendentry{CR}
 \addplot table [x=h, y=eU0TS] {Resu/StokesSinusRT-P1P0-nu6.txt};
 \addlegendentry{Post}
 \end{axis}
\end{tikzpicture}
\begin{tikzpicture}[scale=0.6]
 \begin{axis}[
 height = 8.2cm,
 width = 8.2cm,
 xlabel = {$h$},
 ylabel = {$\eps_0^\nu (p_h)$},
 grid=both,
 major grid style={black!50},
 xmode=log, ymode=log,
 xmin=5e-3, xmax=1e-1,
 ymin=1e-11, ymax=1e-5,
 yticklabel style={
     /pgf/number format/fixed,
     /pgf/number format/precision=0
 },
 scaled y ticks=false,
 xticklabel style={
     /pgf/number format/fixed,
     /pgf/number format/precision=0
 },
 scaled y ticks=false,  
 legend style={at={(1,0)},anchor=south east},
 legend columns=2
 ],
 \addplot table [x=h, y=eP0RT] {Resu/StokesSinusRT-P1P0-nu6.txt};
 \addlegendentry{CR}
 \addplot table [x=h, y=eP0TS] {Resu/StokesSinusRT-P1P0-nu6.txt};
 \addlegendentry{Post}
 \end{axis}      
\end{tikzpicture}
\caption{Sinusoidal velocity, $\nu=10^{-6}$.
Plots of $\eps_0^\nu(\uvec_h)$ (left) and $\eps_0^\nu(p_h)$ (right).}
\label{fig:SinFlowRT6}
\end{figure}

%% file: LowRegFigRT.tex
\begin{figure}[h!]
\centering
\begin{tikzpicture}[scale=0.6]
 \begin{axis}[
 height = 8.2cm,
 width = 8.2cm,
 xlabel = {$h$},
 ylabel = {$\eps_0^\nu(\uvec_h)$},
 grid=both,
 major grid style={black!50},
 xmode=log, ymode=log,
 xmin=5e-3, xmax=1e-1,
 ymin=1e-11, ymax=1e-5,
 yticklabel style={
     /pgf/number format/fixed,
     /pgf/number format/precision=0
 },
 scaled y ticks=false,
 xticklabel style={
     /pgf/number format/fixed,
     /pgf/number format/precision=0
 },
 scaled y ticks=false,  
 legend style={at={(1,0)},anchor=south east},
 legend columns=2
 ],
 \addplot table [x=h, y=eU0RT] {Resu/StokesVortexRT-P1P0-nu6.txt};
 \addlegendentry{CR}
 \addplot table [x=h, y=eU0TS] {Resu/StokesVortexRT-P1P0-nu6.txt};
 \addlegendentry{Post}
 \end{axis}
\end{tikzpicture}
\begin{tikzpicture}[scale=0.6]
 \begin{axis}[
 height = 8.2cm, 
 width = 8.2cm,  
 xlabel = {$h$},
 ylabel = {$\eps_0^\nu(p_h)$},
 grid=both,
 major grid style={black!50},
 xmode=log, ymode=log,
 xmin=5e-3, xmax=1e-1,
 ymin=1e-11, ymax=1e-5,
 yticklabel style={
     /pgf/number format/fixed,
     /pgf/number format/precision=0
 },
 scaled y ticks=false,
 xticklabel style={
     /pgf/number format/fixed,
     /pgf/number format/precision=0
 },
 scaled y ticks=false,  
 legend style={at={(1,0)},anchor=south east},
 legend columns=2
 ],
 \addplot table [x=h, y=eP0RT] {Resu/StokesVortexRT-P1P0-nu6.txt};
 \addlegendentry{CR}
 \addplot table [x=h, y=eP0TS] {Resu/StokesVortexRT-P1P0-nu6.txt};
 \addlegendentry{Post}
 \end{axis}      
\end{tikzpicture}
\caption{Low regularity velocity, $\nu=10^{-6}$.
Plots of $\eps_0^\nu(\uvec_h)$ (left) and $\eps_0^\nu (p_h)$ (right).}
\label{fig:LowRegRT6h}
\end{figure}